\documentclass[a4paper]{article}
\pdfoutput=1
\usepackage{hyperref}
\usepackage{graphicx}
\usepackage{latexsym}
\usepackage{amsfonts}
\usepackage{amssymb}
\usepackage{amsmath}
\usepackage{amsthm}
\usepackage[all,2cell]{xy}
\usepackage{enumerate}
\usepackage{mathrsfs}
\usepackage{verbatim}
\usepackage{textpos}
\usepackage{stmaryrd}
\UseAllTwocells
\input xy

\TPGrid{1}{1} % Set grid to 1 across, 1 down for text-positioning

\theoremstyle{plain}
\newtheorem{thm}{Theorem}
\newtheorem*{thmt}{Theorem}
\newtheorem{cor}[thm]{Corollary}
\newtheorem*{cort}{Corollary}
\newtheorem{lem}[thm]{Lemma}
\newtheorem{prop}[thm]{Proposition}

\theoremstyle{definition}
\newtheorem{defn}[thm]{Definition}

\SelectTips{cm}{}

\newcommand\mi{\text{-1}}
\newcommand{\la}[1]{\stackrel{\,\,\, #1}{\longleftarrow}}
\newcommand{\ra}[1]{\stackrel{\,\,\, #1}{\longrightarrow}}
\newcommand{\sra}[1]{\stackrel{ #1}{\rightarrow}}
\newcommand{\sla}[1]{\stackrel{\,\, #1}{\leftarrow}}
\newcommand\ba{\begin{aligned}}
\newcommand\ea{\end{aligned}}
\newcommand\ig{\includegraphics}

\newcommand\Iso{\stackrel{\sim}{\Rightarrow}}
\newcommand\uIso{\operatorname{uIso}}

\newcommand\id{\operatorname{id}}

\newcommand\ch{\operatorname{ch}}

\newcommand\be{\begin{equation}}
\newcommand\ee{\end{equation}}
\newcommand{\ttr}[1] {{\mathbb{T}r(#1)}}

\newcommand{\Amb}{\operatorname{Amb}}

\newcommand\Sections{\operatorname{Sections}}

\newcommand{\Adj}{\operatorname{Adj}}

\newcommand{\Tr}{\operatorname{Tr}}
\newcommand{\Ind}{\operatorname{Ind}}

\newcommand{\TRep}{\operatorname{2\mathcal{R}ep}}
\newcommand{\LBun}{\operatorname{LBun}}

\newcommand{\Class}{\operatorname{Class}}

\newcommand{\Rep}{\operatorname{Rep}}

\newcommand{\Hom}{\operatorname{Hom}}
\newcommand{\Hilb}{\operatorname{Hilb}}

\newcommand{\Pic}{\operatorname{Pic}}
\newcommand{\Aut}{\operatorname{Aut}}

\newcommand{\End}{\operatorname{End}}
\newcommand{\Vect}{\operatorname{Vect}}
\newcommand{\THilb}{\operatorname{2\mathcal{H}ilb}}
\newcommand{\TVect}{\operatorname{2\mathcal{V}ect}}

\newcommand{\X}{\mathfrak{X}}

\newcommand{\Ob}{\operatorname{Ob}}

\newcommand{\Nat}{\operatorname{Nat}}

\newcommand{\Var}{\operatorname{\mathcal{V}ar}}
\newcommand{\Gerbes}{\operatorname{\mathcal{G}erbes}}

\newcommand{\FFix}{\operatorname{Fix}}
\newcommand{\UOneTor}{U(1)\text{-Tor}}

\newcommand{\arXiv}[1]{\href{http://arxiv.org/abs/#1}{{\tt arXiv:#1}}}

\newcommand{\qalg}[1]{\href{http://arxiv.org/abs/q-alg/#1}{{\tt arXiv:q-alg/#1}}}

\newcommand{\mathQA}[1]{\href{http://arxiv.org/abs/math.QA/#1}{{\tt arXiv:math.QA/#1}}}

\newcommand{\mathCT}[1]{\href{http://arxiv.org/abs/math.CT/#1}{{\tt arXiv:math.CT/#1}}}

\newcommand{\mathDG}[1]{\href{http://arxiv.org/abs/math.DG/#1}{{\tt arXiv:math.DG/#1}}}

\newcommand{\Math}[1]{\href{http://arxiv.org/abs/math/#1}{{\tt arXiv:math/#1}}}

\newcommand{\mathAT}[1]{\href{http://arxiv.org/abs/math.AT/#1}{{\tt arXiv:math.AT/#1}}}

\newcommand{\mathKT}[1]{\href{http://arxiv.org/abs/math.KT/#1}{{\tt arXiv:math.KT/#1}}}

\newcommand{\Fix}[2] {
    \xy
        (0,-0.5)*{#1};
        (0,1)*{\ba \ig{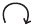} \ea};
        (3,2)*{\scriptstyle #2};
    \endxy}

\newcommand{\sFix}[2] {
    \xy
        (-0.2,-1.7)*{\scriptstyle #1};
        (-0.2,-0.5)*{\ba \ig{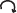}\ea};
        (1.8,0)*{\scriptscriptstyle #2};
    \endxy}

\newcommand{\ssFix}[4] {
    \xy
        (-0.2,-1.7)*{\scriptstyle #1};
        (-0.2,0)*{\ba \ig{sloop.pdf}\ea};
        (#3,#4)*{\scriptscriptstyle #2};
    \endxy}

\newcommand{\Fixx}[4] {
    \xy
        (0,-0.5)*{#1};
        (0,2)*{\ba \ig{loop.pdf}\ea};
        (#3,#4)*{\scriptstyle #2};
    \endxy}

\newcommand{\lag}[2]{ \xy {\ar^(0.45){\scriptscriptstyle #2}_(0.4){\scriptscriptstyle
#1} (6,0)*{}; (0,0)*{}};
\endxy}

 \reversemarginpar
\begin{document}

\title{The geometry of unitary 2-representations of finite groups and their 2-characters}

\author{Bruce Bartlett\\ University of Sheffield\\{\tt \small b.h.bartlett@sheffield.ac.uk}}
\date{}
%\address{University of Sheffield}
%\email{b.h.bartlett@sheffield.ac.uk}
\maketitle

\begin{abstract}
Motivated by topological quantum field theory, we investigate the geometric aspects of unitary 2-representations of finite groups on 2-Hilbert spaces, and their 2-characters. We show how the basic ideas of geometric quantization are `categorified' in this context: just as representations of groups correspond to equivariant line bundles, 2-representations of groups correspond to equivariant gerbes. We also show how the 2-character of a 2-representation can be made functorial with respect to morphisms of 2-representations. Under the geometric correspondence, the 2-character of a 2-representation corresponds to the geometric character of its associated equivariant gerbe. This enables us to show that the complexified 2-character is a unitarily fully faithful functor from the complexified homotopy category of unitary 2-representations to the category of unitary equivariant vector bundles over the group.
\end{abstract}

\section*{Introduction}
One of the main themes of topological quantum field theory (TQFT) is that abstract yet elementary higher categorical ideas can translate into quite sophisticated geometric structures. In this paper we work out a concrete toy model of this sort, in the hope that the same methods will apply in a more advanced setting.

The higher categorical structures we will be concerned with are {\em unitary 2-representations of finite groups}, the {\em 2-category} which they constitute, and their {\em 2-characters}. A 2-representation of a group is a group acting coherently on some sort of linear category. They often arise in practice when a group acts on a geometric or algebraic structure, for then the group will act on the {\em category of representations} of that structure. A 2-representation can be thought of higher-categorically as weak 2-functor from the group (thought of as a one-object 2-category with only identity 2-morphisms) to the 2-category of linear categories. Hence they naturally form a 2-category, if we define the morphisms to be transformations between the weak 2-functors and the 2-morphisms to be modifications (see \cite{ref:leinster_basic_bicategories} for our conventions). Our main point in this paper is that all these higher-categorical ideas translate into concrete geometric structures --- but to understand this correspondence, one first needs to understand the geometric correspondence for {\em ordinary} representations of groups.

\subsubsection*{Geometry of ordinary representations of groups and their characters}
The basic idea of geometric quantization in the equivariant context is that every representation of a group $G$ arises as the `quantization' of a classical geometric system having symmetry group $G$. Normally this is expressed in the language of symplectic geometry and polarizations \cite{ref:woodhouse}, but one can find an elementary categorical formulation of it in a simple setting which for our purposes still displays many of the essential features, as we now explain (we apologize to the reader for being somewhat sketchy in this subsection; an explicit write-up will appear elsewhere \cite{ref:bartlett_ord}. We remark that our approach is in the spirit of the {\em coherent states} formalism as in \cite{ref:kirwin, ref:spera}).

\begin{figure}
\centering
\ig{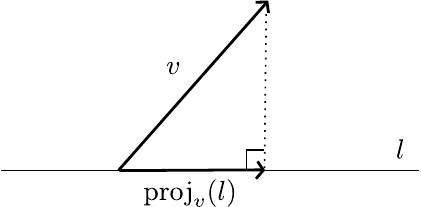}
\caption{\label{projfig}The section of the line bundle over projective space associated to a vector $v \in V$.}
\end{figure}

The main thing to understand is that {\em every finite dimensional Hilbert space $V$ identifies antilinearly as the space of sections of a holomorphic line bundle}. Indeed, given $V$ we have the associated line bundle over its projective space $\tau_V \rightarrow \mathbb{P}(V)$ whose fiber at a line $l \in \mathbb{P}(V)$ is the line $l$ itself. Thus to a vector $v \in V$ we may assign a holomorphic section of $\tau_V$ by orthogonally projecting $v$ onto every line $l$ (see Figure \ref{projfig})\footnote{In more orthodox terminology we are using the inner products to identify $\tau_V$ with the dual of the tautological line bundle.}, and all holomorphic sections of $\tau_V$ are of this form. A category theorist should think of this as the `decategorified Yoneda lemma' since it says that a vector is determined by all of its inner products.

One can make this correspondence equivariant and also upgrade it to an equivalence of categories in the following way. Given a compact Lie group $G$, define the category $\LBun(G)$ as follows. An object is a holomorphic equivariant unitary hermitian line bundle $L \rightarrow X$ over a compact hermitian manifold $X$ acted on isometrically by $G$. That is, $L$ is a holomorphic collection of hermitian lines $\{L_x\}$ where $x$ ranges over $X$, equipped with unitary maps
 \[
   L_x \rightarrow L_{g \cdot x}
 \]
for each $g \in G$ lifting the action of $G$ on $X$. A morphism from $L \rightarrow X$ to $Q \rightarrow Y$ is an {\em equivariant kernel} --- an equivariant holomorphic collection of linear maps
 \[
 \langle y | E | x \rangle \colon L_x \rightarrow Q_y
 \]
from each fiber of $L$ to each fiber of $Q$, and one composes these by integration --- that is, if $E'$ is another equivariant kernel from $Q \rightarrow Y$ to $R \rightarrow Z$ then
  \be \label{kerncomp}
   \langle z | E' \circ E | x \rangle = \int_Y dy \; \langle z | E' | y \rangle \circ \langle y | E | x \rangle.
  \ee
In this setting the idea of geometric quantization translates into the statement that there is an equivalence between the category of unitary representations of $G$ and the category of equivariant line bundles:
 \begin{align*}
 \Rep(G) & \rightarrow \LBun(G)\\
 V & \mapsto \tau_V \\
 \overline{\Gamma(L)}  & \mapsfrom L
 \end{align*}
There is also a notion of the {\em geometric character} $\ch_L$ of an equivariant line bundle $L$, defined by integrating the group action over the manifold. One first repackages the equivariant line bundle as a kernel
 \[
   \langle y | g | x \rangle \colon L_x \rightarrow L_y
 \]
and then one defines the geometric character of $g \in G$ as the integral
 \[
  \ch_L (g) := \int_X dx \, \Tr \langle x | g | x \rangle.
 \]
Under strong enough conditions on $L$ and $X$ this integral localizes over the {\em fixed points} of $g$ on $X$ (for projective space, these are precisely the {\em eigenlines}, so one is simply summing the eigenvalues) --- a statement we will make more precise in \cite{ref:bartlett_ord}.

In any event, it is not hard to show that that the character of a representation corresponds to the geometric character of its associated equivariant line bundle; in other words, we have the commutative diagram
 \[
  \xymatrix @1 @C=0.1in { [\Rep(G)]_\mathbb{C} \ar[rr]^\cong \ar[dr]_{\chi}  && [\LBun(G)]_\mathbb{C} \ar[dl]^{\ch} \\ & \Class(G)}
 \]
where $[ \cdot ]_\mathbb{C}$ refers to the complexified Grothendieck groups of these categories and $\Class(G)$ is the space of class functions on the group. Moreover, all the maps above are {\em unitary} with respect to the natural inner products involved.

\subsubsection*{Categorifying the geometric correspondence}
The main point of this paper is to show that the above geometric correspondence `categorifies' appropriately to the setting of unitary {\em 2}-representations. By a `unitary 2-representation' we will mean a group acting unitarily and coherently on a {\em 2-Hilbert space} (see \cite{ref:baez_2_hilbert_spaces}). These are abelian linear categories equipped with a duality and a compatible inner-product; they stand in the same relation with Kapranov and Voevodsky's {\em 2-vector spaces} \cite{ref:kapranov_voevodsky} as finite-dimensional Hilbert spaces are to finite-dimensional vector spaces. The main thing to bear in mind is that these categories are {\em semisimple}, and this means that the geometric structures we will derive in this paper will always be discrete --- in particular, we are obliged to restrict ourselves to 2-representations of {\em finite} groups. Nevertheless, we hope that similar ideas will apply in the non-semisimple context, such as group actions on derived categories of sheaves \cite{ref:deligne, ref:ganter_kapranov_rep_char_theory}.

\subsubsection*{2-characters}
The first thing to do is to develop a good theory of 2-characters. Just as the character of an ordinary representation of a group is an assignment of a {\em number} to each element of the group, invariant under conjugation, the {\em 2-character} of a 2-representation is an assignment
of a {\em vector space} to each element of the group, together with specified isomorphisms relating the vector spaces assigned to conjugate group elements.

Now 2-characters were introduced independently by Ganter and Kapranov \cite{ref:ganter_kapranov_rep_char_theory} while we were working on this project; however we develop a number of results about them not present in \cite{ref:ganter_kapranov_rep_char_theory}. Our first step of departure is to use {\em string diagrams} as a convenient notation for working with 2-representations and their 2-characters; we explain this notation in the opening section. Using this notation, we show that not only can one take the 2-character of a 2-representation, but one can also take the 2-character of a {\em morphism} of 2-representations. To do this, one needs to have good control over the ambidextrous adjunctions in the 2-category of 2-Hilbert spaces, or geometrically speaking, one needs to ensure that one can choose a flat section of the `ambidextrous adjunction bundle'. We call this an {\em even-handed structure}, and the behaviour of the 2-character on morphisms will in general {\em depend} on the choice of this structure. However, we show that the 2-category of 2-Hilbert spaces has a canonical even-handed structure, which uses the inner products and duality on the hom-sets in an essential way (this is analogous to the way that the adjoint of a linear map between vector spaces requires an inner product). This is one of the features of working with unitary 2-representations which is not available for 2-representations on unadorned 2-vector spaces.

Just as the ordinary character of a representation does not depend on the isomorphism class of the representation, the 2-character of a {\em morphism} of 2-representations does not depend on its isomorphism class, hence it descends to a functor from the {\em homotopy category} of unitary 2-representations to the category of equivariant vector bundles over the group:
   \[
    \chi \colon [\TRep(G)] \rightarrow \Hilb_G (G).
    \]
Our main result in this paper is that after one tensors the hom-sets in $[\TRep(G)]$ with $\mathbb{C}$, the resulting functor is {\em unitarily fully faithful}. This is the categorification of the fact that the ordinary character is a unitary map from the complexified Grothendieck group of unitary representations to the space of class functions. To prove this result, we will develop a geometric correspondence for unitary 2-representations in terms of {\em finite equivariant gerbes} analogous to the correspondence between ordinary representations and equivariant line bundles above. Then we show that under this correspondence, the 2-character of the 2-representation corresponds to the {\em geometric character} of the associated equivariant gerbe, from which the result follows since we can use a theorem of Willerton \cite{ref:simon}, developed in the context of twisted representations of finite groupoids, to obtain a detailed understanding of geometric characters.

\subsubsection*{Equivariant gerbes}
For our purposes, a {\em finite equivariant gerbe} is a $U(1)$-central extension $\X$ of the action groupoid $X_G$ associated to a finite $G$-set $X$ (see \cite{ref:brylinski, ref:behrend_xu, ref:simon} for background material on gerbes). That is, $\X$ is a groupoid having the same objects as $X_G$ with the property that each morphism in $X_G$ is replaced by a $U(1)$-torsor worth of morphisms in $\X$. An equivariant gerbe can also be thought of as representing a `equivariant hermitian 2-line bundle' over $X$; one way to see this is that by choosing a set-theoretic section of the gerbe one can extract a $U(1)$-valued groupoid 2-cocycle $\phi$, which can be used to form a fibration of categories $X_G \times_\phi \Hilb \rightarrow X_G$ (see \cite{ref:simon, ref:brylinski_mclaughlin, ref:cegarra_et_al_graded_extensions}). Since we will be dealing with integration in various forms, we will also equip our gerbes with {\em metrics}. Equivariant gerbes form a 2-category $\Gerbes(G)$, if we define a 1-morphism $\X \rightarrow \X'$ between gerbes to be a unitary equivariant vector bundle over their product $\X' \otimes \overline{\X}$ (this should be thought of as a {\em categorified kernel}) and a 2-morphism to be a morphism of equivariant vector bundles.

\subsubsection*{Unitary 2-representations and equivariant gerbes}
It turns out that from a marked unitary 2-representation of $G$ (a 2-Hilbert space is {\em marked} if it is endowed with distinguished simple objects), one can extract a finite equivariant gerbe equipped with a metric, and similarly for the morphisms and 2-morphisms, leading to an equivalence of 2-categories
 \[
  \TRep_m (G) \stackrel{\sim}{\rightarrow} \Gerbes(G).
 \]
This should be regarded as a `categorification' (at least in our finite discrete setting) of the aforementioned equivalence between the category of unitary representations of $G$ and the category of equivariant line bundles.

\subsubsection*{Geometric characters of equivariant gerbes}
We have seen that the geometric character of an equivariant line bundle can often be expressed as a sum over the fixed points of the group action. Similarly we define the geometric character of an equivariant gerbe as the {\em space of sections} over the fixed points of the gerbe (readers familiar with these ideas will recognize this as the push-forward of the transgression map in our context, as in \cite{ref:tu_xu_laurent_gengoux}). Since the group acts on these fixed points by conjugation, the geometric character also produces an equivariant vector bundle over the group. We show how one can also apply the geometric character to a {\em morphism} of equivariant gerbes, so as to obtain a morphism between the corresponding equivariant vector bundles. Just as for 2-representations, the geometric character on morphisms only depends on their {\em isomorphism classes}, hence it also descends to a functor from the homotopy category of equivariant gerbes to the category of equivariant vector bundles over the group,
 \[
  \ch \colon [\Gerbes(G)] \rightarrow \Hilb_G(G).
 \]
Now, the morphisms in $[\Gerbes(G)]$ may be regarded as twisted representations of groupoids, and it turns out that at the level of morphisms the geometric character functor essentially takes the {\em twisted characters} of these representations. Thus one can apply the technology of Willerton \cite{ref:simon} to conclude that after one tensors the hom-sets in $[\Gerbes(G)]$ with $\mathbb{C}$, the geometric character functor is {\em unitarily fully faithful}.

\subsubsection*{Main result}
As we have explained, our main result is a `categorification' of the aforementioned geometric correspondence for ordinary representations.
\begin{thmt} The 2-character of a marked unitary 2-representation is unitarily naturally
isomorphic to the geometric character (i.e. the push-forward of the transgression) of its associated equivariant
gerbe:
 \[
 \xymatrix @1 @C=0.1in{[\TRep_m(G)] \ar[dr]_{\chi} \ar[rr]^{\simeq} && [\Gerbes(G)] \ar[dl]^{\ch} \\ & \Hilb_G(G)}
 \]
\begin{textblock}{0.3}(0.245,-0.0515)
  $   \ba \ig{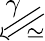} \ea $
\end{textblock}
\end{thmt}
\begin{cort} The complexified 2-character functor
 \[
  \chi_\mathbb{C} \colon [\TRep(G)]_\mathbb{C} \rightarrow \Hilb_G(G)
  \]
is a unitarily fully faithful functor from the complexified homotopy category of unitary 2-representations of $G$ to the category of unitary equivariant vector bundles over $G$.
\end{cort}

\subsubsection*{Comparison with previous work} There have already been a number of works on 2-representations (see eg. \cite{ref:barrett_mackaay, ref:crane_yetter, ref:deligne, ref:elgueta, ref:ganter_kapranov_rep_char_theory, ref:morton, ref:ostrik}), the most relevant for this paper being that of Elgueta \cite{ref:elgueta} and that of Ganter and
Kapranov \cite{ref:ganter_kapranov_rep_char_theory}. Elgueta performed a thorough
and careful investigation of the 2-category of 2-representations of a 2-group (a 2-group is a monoidal category with structure and properties analogous to that of a group \cite{ref:baez_lauda_2-groups}) acting on Kapranov and Voevodsky's 2-vector spaces, and his motivation was therefore to work with co-ordinatized versions of 2-vector spaces amenable for direct computation, and to classify the various structures which appear.

Ganter and Kapranov were motivated by equivariant homotopy theory, namely to try and find a categorical construction which would produce the sort of generalized group characters which crop up in Morava $E$-theory; they (independently) introduced the {\em categorical character} (which we call the 2-character) of a 2-representation and showed that it indeed achieves this purpose. Since they had no reason not to, they also worked with co-ordinatized 2-vector spaces (this time of the form $\Vect^n$); also they did not investigate in any depth morphisms and 2-morphisms of 2-representations. This present paper was in preparation at the time \cite{ref:ganter_kapranov_rep_char_theory} appeared and the author apologizes for the lengthy delay.

Our motivation has been extended topological quantum field theory, where the 2-category of 2-representations of a group appears as the `2-category associated to the point' in a finite version of Chern-Simons theory called the {\em untwisted finite group model} (see \cite{ref:dijkgraaf_witten, ref:freed, ref:morton, ref:simon, ref:baez_dolan_HDA0} for background). We remark here that in the {\em twisted} model (the twisting is given by a group 3-cocycle $\omega \in Z^3(G, U(1)$) the `2-category assigned to the point' is the 2-category of 2-representations of the {\em 2}-group $G_\omega$ constructed from $G$ and $\omega$. In this paper we have restricted ourself to the untwisted case, since the geometry of the twisted model is a bit more intricate (one must essentially replace equivariant gerbes by {\em twisted} equivariant gerbes), and also because the untwisted setting makes for a cleaner analogy between the geometric picture of ordinary unitary representations and that of unitary 2-representations. We hope to study more carefully the geometry of the twisted model in future work; for now we refer the reader to the recent \cite{ref:baez_baratin_freidel_wise} which studies strict 2-representations of Lie 2-groups on `higher Hilbert spaces' (these are categories whose objects are `measurable fields of Hilbert spaces' supported over a measurable space $X$, a setting which allows for continuous geometry).

In any event, the language of Chern-Simons theory is the geometric language of moduli-stacks, line-bundles, equivariant structures, flat sections and such like, and this has therefore motivated our
approach to 2-representations and is what distinguishes our approach from previous approaches (though we remark that related ideas do appear in \cite{ref:ganter_kapranov_rep_char_theory}).   For instance, as far as possible we try to work directly with the underlying 2-Hilbert spaces of the 2-representations themselves as opposed to some `co-ordinitization' of them, a strategy which might be important in a more intricate geometric setting.  We hope that some of the ideas we have developed in this paper will also translate into the more advanced geometric contexts of \cite{ref:ganter_kapranov_rep_char_theory}, as well as \cite{ref:freed_teleman_hopkins}.

\subsubsection*{Overview of paper}
In Section \ref{stsec} we remind the reader of how the string diagram notation for 2-categories works. In Section \ref{utrep} we recall the notion of a 2-Hilbert space due to Baez \cite{ref:baez_2_hilbert_spaces}, and we define what we mean by the 2-category of unitary 2-representations, expanding out all the definitions in terms of string diagrams, where they take a particularly simple form. We also give a number of examples of unitary 2-representations. Finally we develop the idea of an {\em even-handed structure} on a 2-category as a consistent choice of isomorphism classes of ambidextrous adjoints, and we show that the 2-category of 2-Hilbert spaces has a canonical such structure. We also explain what it means for a 2-representation to be compatible with a given even-handed structure.

In Section \ref{charsecc} we use the string diagram technology to define 2-characters of unitary 2-representations, and we show how to make this construction functorial with respect to morphisms of 2-representations.

In Section \ref{gggsec} we introduce the main geometric notions of this paper: finite equivariant gerbes equipped with metrics, the 2-category which they constitute, and the twisted character of an equivariant vector bundle over a gerbe. In Section \ref{gch} we use this language to define the geometric character of an equivariant gerbe, and we show how to make this construction functorial with respect to morphisms of equivariant gerbes.

In Section \ref{togsec} we explain how to extract an equivariant gerbe from a unitary 2-representation, and similarly for the morphisms and 2-morphisms, leading to an equivalence of 2-categories. Finally in Section \ref{2char2rep} we bring together all these concepts, and show that the 2-character of a 2-representation corresponds naturally to the geometric character of its associated equivariant gerbe. We use this to conclude that the 2-character is a unitarily fully faithful functor after one passes to the complexified homotopy category.

\subsubsection*{Acknowledgements} During the preparation of this paper I have benefited from discussions with and valuable remarks from Urs Schreiber, James Cranch, John Baez, Jeffrey Morton, Nora Ganter, Mikhail Kapranov, David Gepner, Richard Hepworth, Eugenia Cheng, Tom Bridgeland, Vic Snaith, Ieke Moerdijk, Frank Neumann and Mathieu Anel. A special acknowledgement goes to my supervisor Simon Willerton, who first suggested to me the topic of 2-representations and introduced me to the notion of the 2-character; much of this present work can be seen as a follow-up to \cite{ref:simon}.  I also wish to acknowledge support from the Excellence Exchange Scheme at the University of Sheffield. Finally I wish to thank the organizers of the Max Kelly conference in Cape Town for a special and memorable occasion.

\section{String diagrams\label{stsec}} In this section we briefly recall the string diagram notation for 2-categories. This notation is particularly suited to describe structures such as adjunctions and monads, and we will find it very useful when we discuss 2-characters in Section \ref{charsecc}.

\begin{figure}
 \centering
 \ig{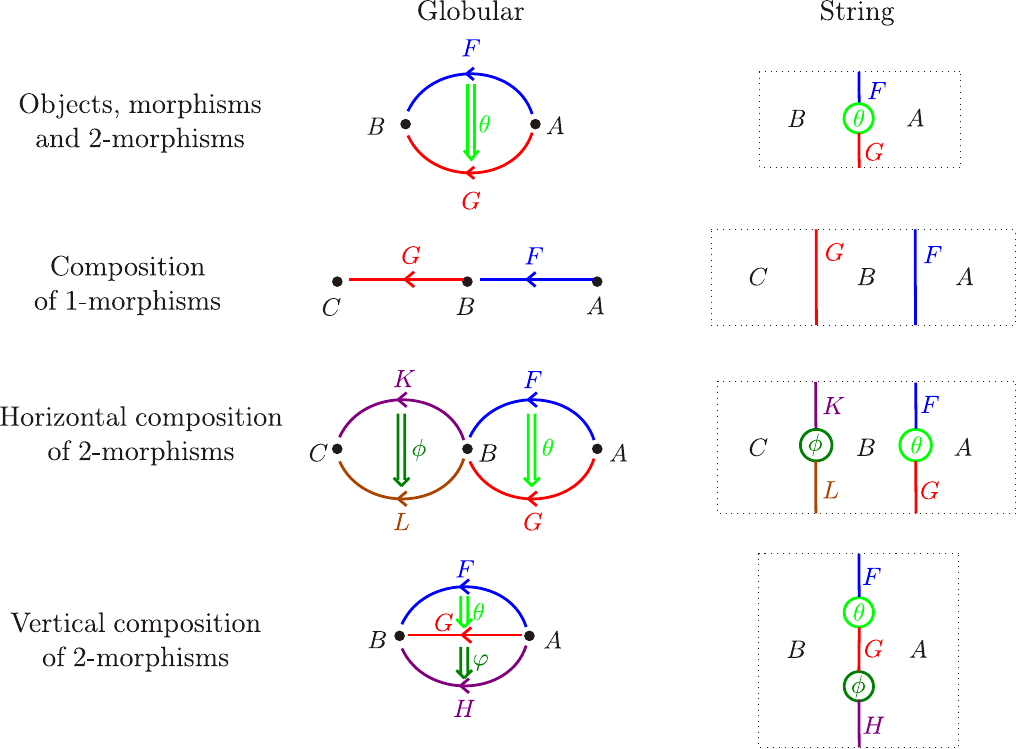}
 \caption{\label{fig:sd} Globular notation for 2-categories versus string diagram notation. }
 \end{figure}

String diagrams are a two-dimensional graphical notation for working with 2-categories, and may be regarded as the `Poincar\'{e} duals' of the ordinary globular notation. The basic idea is summarized in Figure \ref{fig:sd}. The reader who is still confused by these diagrams is referred to Section 2.2 of \cite{ref:lauda}, Section 1.1 of \cite{ref:caldararu_willerton} or \cite{ref:joyal_street, ref:street_string_diagrams} for more details. In our diagrams, composition of 1-morphisms runs from right to left (so a $G$ to the left of an $F$ means $G$ after $F$), and composition of 2-morphisms runs from top to bottom. We stress that string diagrams are not merely a mnemonic but are a perfectly rigorous notation.

We make the important remark here that by a `2-category' we mean a {\em not-necessarily-strict} 2-category (also called a {\em bicategory}). In this paper we will almost exclusively be using string diagrams in the context of the 2-category of 2-Hilbert spaces, which {\em is} a strict 2-category. However we will have occasion to use them in the context of a weak 2-category when we discuss even-handed structures for general 2-categories in Section \ref{EvenSec}. Let us therefore make some remarks on the interpretation of these diagrams when the 2-category is not strict.

A string diagram is a graphical notation for a specific 2-morphism in a 2-category. The source and target 1-morphisms of this 2-morphism are the composites of the top and bottom 1-morphisms represented in the diagram, respectively. When the 2-category is weak, one therefore needs to be given the additional information of precisely how these composites which make up the source and target 1-morphisms are to be parenthesized (this might include arbitrary insertions of identity 1-morphisms). However, once a parenthesis choice has been made for the source and target 1-morphisms, coherence for 2-categories --- in the form which says `all diagrams of constraints commute' as in Chapter 1 of the thesis of Gurski \cite{ref:Gurski} --- implies that the resulting 2-morphism represented by the diagram is {\em unique}, and does not depend on the choice of parentheses, associators and unit 2-isomorphisms used to interpret the interior of the diagram.

Now, whenever a string diagrams occurs in this paper it will always be manifestly clear what the input and output 1-morphisms are. For instance:
\[
\text{``The 2-morphism $\eta\colon \id \Rightarrow (G \circ F) \circ (F^* \circ G^*)$ is defined as } \ba \ig{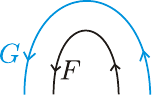} \ea \text{ .''}
   \]
In other words, every string diagram in this paper has a precise and rigorous meaning.

\section{Unitary 2-representations\label{utrep}} We begin this section by recalling the notion of a 2-Hilbert space due to Baez \cite{ref:baez_2_hilbert_spaces}. Then we define the 2-category of unitary 2-representations of a finite group. We do this first in an abstract higher-categorical way, and then we spell out this definition in terms of string diagrams by introducing graphical elements to depict the various pieces of data involved. We give some examples of 2-representations, such as those arising from exact sequences of groups. Finally we show that $\THilb$ comes equipped with a canonical {\em even-handed structure} --- a coherent way to equip every right adjoint also as a left adjoint. This will be necessary in the next section in order to show that the 2-character can be made functorial.

\subsection{2-Hilbert spaces}
A {\em 2-Hilbert space} is a `categorification' of a Hilbert space; that is, it is a category with structure and properties analogous to those of a Hilbert space. They were introduced by Baez \cite{ref:baez_2_hilbert_spaces}, and we need them because we want to work with {\em unitary}
2-representations. They are similar to Kapranov and Voevodsky's {\em 2-vector spaces}, the
difference between them being analogous to the difference between a Hilbert space and an ordinary vector space. We remark that by the words `Hilbert space' we will inevitably mean simply a finite-dimensional inner product space; this terminology is standard in a topological quantum field theory context \cite{ref:freed}.

A good example of a 2-Hilbert space is the category of unitary finite-dimensional representations of a finite group (more generally, one could consider twisted unitary representations of a finite group{\em oid}). Observe that this category is especially easy to understand: every representation is a direct sum of irreducible representations. Although it is equivalent to the category of finite-dimensional representations with no inner products involved (which is a 2-vector space in the sense of Kapranov and Voevodsky), it has extra structure because every morphism $f \colon V \rightarrow W$ has an adjoint $f^* \colon W \rightarrow V$ and the hom-sets have a natural inner product, $(f,g) = \frac{1}{|G|}\Tr(f^*g)$. This illustrates the difference between a 2-Hilbert space and a 2-vector space.

\subsubsection{The definition}
We write Hilb for the category of finite dimensional Hilbert spaces and linear maps. A Hilb-category is a category enriched over Hilb.
\begin{defn} A {\em 2-Hilbert space} is an abelian Hilb-category $H$, equipped with antilinear maps $* \colon
 \text{hom}(x,y) \rightarrow \text{hom}(y,x)$ for all $x, y \in H$,
 such that
  \begin{itemize}
   \item $f^{* *} = f$,
   \item $(f g)^* = g^* f^*$,
   \item $(f g, h) = (g, f^* h)$,
   \item $(f g, h) = (f, h g^*)$
   \end{itemize}
  whenever both sides of the equation are defined.
\end{defn}
The thing to keep in mind about 2-Hilbert spaces is that they are {\em automatically semisimple} --- that is, there exist objects $e_i$ labeled by a set $I$ such that $\Hom(e_i, e_j) \cong \delta_{ij}k$ (such objects are called {\em simple}) and such that for any two objects $x$ and $y$ the composition map
  \[
   \Hom(x,y) \leftarrow \bigoplus_{i \in I} \Hom(x, e_i) \otimes \Hom(e_i, y).
  \]
is an isomorphism.
\begin{prop}[{Baez \cite{ref:baez_2_hilbert_spaces}}] Every 2-Hilbert space is semisimple.
\end{prop}

The reason is that the arrow algebra associated to a finite set of objects in a 2-Hilbert space forms an $H^*$-algebra (a finite-dimensional algebra with an inner product and a compatible antilinear involution), and it is a result of Ambrose \cite{ref:ambrose} that such algebras are weighted direct sums of matrix algebras. We see that {\em it
is precisely the geometric ingredient of `duality'} (the inner products and the $*$-structure) which causes
2-Hilbert spaces to be semisimple.

This is the main conceptual difference between 2-Hilbert spaces and 2-vector spaces. The latter were defined by Kapranov and Voevodsky essentially as `an abelian $\Vect$-module category equivalent to $\Vect^n$ for some $n \in \mathbb{N}$'. A
good feature of their definition is that it explicitly includes a prescription for categorified scalar
multiplication, an ingredient which is missing from the definition of 2-Hilbert spaces, but which can be useful
for some constructions (on the other hand, it can easily be added in). A disappointing feature of their
definition though is that it adds in semisimplicity by hand, whereas 2-Hilbert spaces are defined {\em
intrinsically}, and semisimplicity is a consequence.

For instance, the category $\Rep(A)$ of finite-dimensional representations of an algebra $A$ is always an
abelian category with a $\Vect$-module structure, but it fails to be semisimple in general, precisely because of
the lack of duality in $\Rep(A)$.

\subsubsection{The 2-category of 2-Hilbert spaces\label{hilb2}}
A 2-Hilbert space $H$ is called  {\em finite dimensional} if there are only a finite number of non-isomorphic simple
objects; this number is called the {\em dimension} of $H$. A functor $F \colon H \rightarrow H'$ between 2-Hilbert spaces is called {\em linear} if it is linear on
the level of hom-sets and preserves direct sums, in the sense that if $x \oplus y$ is a direct sum of $x, y \in
H$, then $F(x \oplus y)$ is a direct sum of $F(x), F(y) \in H'$. It is called a {\em $*$-functor} if $F(f^*) = F(f)^*$ for all morphisms $f$ in $H$.

\begin{defn} The 2-category $\THilb$ of 2-Hilbert spaces has finite-dimensional 2-Hilbert spaces for objects,
linear $*$-functors for morphisms, and natural transformations for 2-morphisms.
\end{defn}

A natural transformation $\theta \colon F \Rightarrow F'$ between morphisms $F, F' \colon H \rightarrow H'$ of 2-Hilbert spaces is called {\em unitary} if all its components are unitary --- that is, if $\theta_x^* \theta_x = \id_{F(x)}$ and $\theta_x \theta_x^* = \id_{F'(x)}$ for all $x \in H$. A pair $H, H'$ of 2-Hilbert spaces are called {\em unitarily equivalent} if there are
linear $*$-preserving functors $F \colon H \rightarrow H'$ and $G \colon H' \rightarrow H$ together with unitary natural
isomorphisms $\eta \colon \id_H \Iso GF$ and $\epsilon \colon FG \Iso \id_{H'}$ forming an adjunction. We will call them {\em strongly} unitarily equivalent (this refined notion is not considered in \cite{ref:baez_2_hilbert_spaces}) if the functors $F$ and $G$ are also unitary linear maps at the level of hom-sets.

The important things to remember about 2-Hilbert spaces, and the morphisms and 2-morphisms between them, are the following:
 \begin{itemize}
 \item A 2-Hilbert space $H$ is determined up to unitary equivalence simply by its dimension (see \cite{ref:baez_2_hilbert_spaces}). It is determined up to {\em strong} unitary equivalence by its dimension and the {\em scale factors} on the simple objects --- the positive real numbers $k_i = (\id_{e_i}, \id_{e_i})$ (if $H$ is the category of unitary representations of a finite group, these numbers are the dimensions of the irreducible representations divided by the order of the group). These scale factors are the extra information not present in an unadorned semisimple category.

 \item Since 2-Hilbert spaces are semisimple, a linear $*$-functor $F \colon H \rightarrow H'$ between them is determined up to unitary natural isomorphism simply by the vector spaces $\Hom(e_\mu, Fe_i)$, where $e_\mu$ and $e_i$ run over a choice of simple objects for $H'$ and $H$ respectively.

 \item Similarly, a natural transformation $\theta \colon F \Rightarrow G$ is freely and uniquely determined by its components $\theta_{e_i} \colon Fe_i \rightarrow Ge_i$ on the simple objects $e_i$. This gives the vector space of natural transformations an inner product via the formula (compare equation (3.16) of \cite{ref:freed2})
      \[
       \langle \theta,  \theta' \rangle = \sum_i k_i (\theta_{e_i}, \theta'_{e_i}).
      \]
 \end{itemize}
In this way, one can view $\THilb$ as a discrete version of the 2-category $\Var$ (see \cite{ref:caldararu_willerton, ref:ganter_kapranov_rep_char_theory}), whose objects are the derived categories $D(X)$ of coherent sheaves over smooth projective algebraic varieties $X$ and whose hom-categories are the derived categories over the product $Y \times X$, 
 \[
  \Hom_{\Var} (X, Y) = D (Y \times X).
 \]

\subsubsection*{Decategorifying a 2-Hilbert space}
If $H$ is a 2-Hilbert space, we will write $[H]_\mathbb{C}$ for the {\em complexified Grothendieck group} of $H$ --- the tensor product of $\mathbb{C}$ with the abelian semigroup generated by the isomorphism classes of objects $[v]$ in $H$ under the relations $[v \oplus w] = [v] + [w]$. A basis for $[H]_\mathbb{C}$ is given by the isomorphism classes $[e_i]$ of simple objects of $H$. We regard $[H]_\mathbb{C}$ as a Hilbert space with inner product defined on the generating elements by  
 \[
 ([v], [w]) = \dim \Hom(v, w).
\]

\subsection{2-representations in terms of string diagrams}
We now define the 2-category $\TRep(G)$ of unitary 2-representations of a finite group $G$. We do this in two stages --- firstly we define it in a terse higher-categorical way, and then we expand out this definition explicitly in traditional notation as well as in string diagrams, introducing new graphical elements to depict the various pieces of data involved.

Since there are various conventions for terminology for 2-categories, we remark that we are essentially using those of Leinster \cite{ref:leinster_basic_bicategories}. The reader is assured that those parts of the definition below mentioning the word `unitary' will be explained shortly.

\begin{defn}[{compare \cite{ref:elgueta, ref:crane_yetter, ref:barrett_mackaay,  ref:ganter_kapranov_rep_char_theory, ref:ostrik}}] The 2-category $\TRep(G)$ of unitary 2-representations of a finite group $G$ is defined as follows. An object is a unitary weak 2-functor $G \rightarrow \THilb$ (where $G$ is thought of as a 2-category which has only one object, with the elements of $G$ as 1-morphisms, and only identity 2-morphisms). A morphism is a
transformation whose coherence isomorphisms are unitary, and a 2-morphism is a modification.
\end{defn}

We now expand this definition out.
%XXX check labels for equations, which get used and which don't.
\subsubsection{Unitary 2-representations}
A unitary 2-representation of $G$ consists of:
\begin{itemize}
   \item A finite-dimensional 2-Hilbert space $H$, drawn as
    \[
     \ba \ig{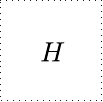}
   \ea \quad \text{or} \quad \ba \ig{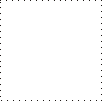} \ea \, \, \text{\; when $H$ is
   understood,}
   \]
   \item For each $g \in G$, a linear $*$-functor $H \la{\alpha_g} H$ which at the level of hom-sets is a unitary linear map, drawn as
   \[
     \ba \ig{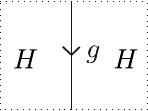} \ea \text{\; or simply \;} \ba \ig{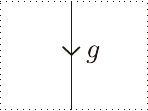} \ea,
    \]
   \item A unitary 2-isomorphism $\phi(e) \colon \id_H \Rightarrow \alpha_e$ (where $e$ is the identity element of $G$), and for each $g_1, g_2 \in G$, a unitary 2-isomorphism $\phi(g_2,g_1)
   \colon \alpha_{g_2} \circ \alpha_{g_1} \Rightarrow \alpha_{g_2g_1}$, drawn as
    \[
     \ba
     \ig{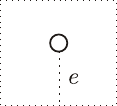} \ea \quad \ba \xymatrix @R=0.5cm {\id_A \ar@{=>}[d]^{\phi(e)}
     \\ \alpha_g} \ea \qquad , \qquad
    \ba \ig{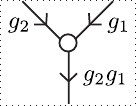} \ea \quad \ba \xymatrix @R=0.5cm {\alpha_{g_2} \circ
     \alpha_{g_1} \ar@{=>}[d]^{\phi(g_2,g_1)} \\ \alpha_{g_2 g_1}} \ea
     \]
\end{itemize}
such that
 \[
  \ba \xymatrix{ & \alpha_g \ar@{=>}[dl]_{\phi(e) \ast \id}
  \ar@{=>}[dr]^{\id \ast \phi(e)} \ar@{=}[dd] \\ \alpha_e \circ \alpha_g \ar@{=>}[dr]_{\phi(e,g)} & &
  \alpha_g \circ \alpha_e \ar@{=>}[dl]^{\phi(g,e)} \\ & \alpha_g} \ea \quad
  \text{and} \quad \ba
  \xymatrix{ & \alpha_{g_3} \circ \alpha_{g_2} \circ \alpha_{g_1} \ar@{=>}[dr]^{\phi(g_3,g_2) \ast \id} \ar@{=>}[dl]_{\id \ast \phi(g_2,g_1)} \\ \alpha_{g_3 g_2} \circ \alpha_{g_1} \ar@{=>}[dr]_{\phi(g_3, g_2 g_1)} &
  & \alpha_{g_3} \circ \alpha_{g_2 g_1} \ar@{=>}[dl]^{\phi(g_3 g_2, g_1)} \\ & \alpha_{g_3 g_2 g_1}} \ea
 \]
commute, or in string diagrams,
  \be
  \label{br2} \ba \ig{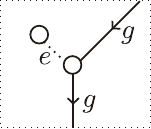} \ea = \ba \ig{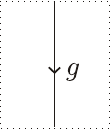} \ea = \ba \ig{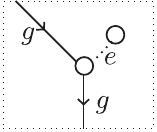} \ea \qquad \text{and} \qquad \ba \ig{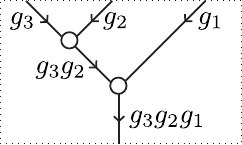} \ea = \ba \ig{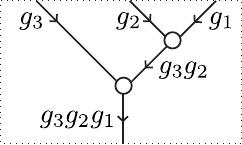} \ea.
  \ee
We will draw the inverse 2-isomorphisms $\phi(e)^* \colon \alpha_e \Rightarrow \id_H$
and $\phi(g_2,g_1)^*  \colon \alpha_{g_2 g_1} \Rightarrow \alpha_{g_2} \circ \alpha_{g_1}$ as
 \[
  \ba \xymatrix @R=0.5cm {\alpha_e \ar@{=>}[d]^{\phi(e)^*} \\ \id_A} \ea \quad \ba \ig{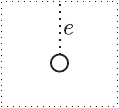} \ea \qquad , \qquad \ba \xymatrix @R=0.5cm {\alpha_{g_2 g_1}
  \ar@{=>}[d]^{\phi(g_2,g_1)^*} \\ \alpha_{g_2} \circ \alpha_{g_1}} \ea \quad \ba \ig{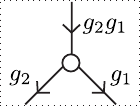} \ea .
 \]
The fact that these satisfy $\phi(e)^* \phi(e) = \id$ and $\phi(e) \phi(e)^* = \id$, and similarly for the $\phi(g_2, g_1)$, is drawn as follows:
 \begin{eqnarray}
  \ba \ig{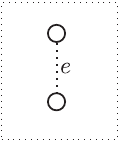} \ea = \ba \ig{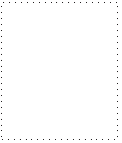} \ea \quad &,& \quad \ba
  \ig{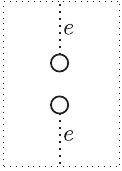} \ea = \ba \ig{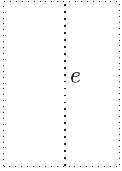} \ea \label{zr1}  \\
  \label{br1} \ba \ig{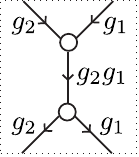} \ea = \ba \ig{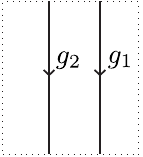} \ea \quad &,& \quad \ba
  \ig{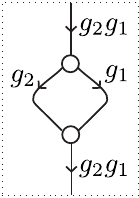} \ea = \ba \ig{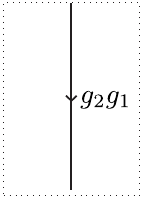} \ea. \label{zr2}
 \end{eqnarray}
We will abbreviate all of this data $(H, \{\alpha_g\}, \phi(e), \{\phi(g_2, g_1)\})$ simply as $\alpha$.

\subsubsection{Morphisms} A morphism $\sigma \colon \alpha \rightarrow \beta$ of unitary 2-representations is a transformation from $\alpha$ to $\beta$ whose coherence isomorphisms are unitary. Thus, if $\alpha = (H_\alpha, \{\alpha_g\}, \phi(e), \{\phi(g_2, g_1)\})$ and $\beta = (H_\beta, \{\beta_g\}, \psi(e), \{\psi(g_2, g_1\})$, then it consists of
 \begin{itemize}
  \item A linear $*$-functor $\sigma \colon H_\alpha \rightarrow H_\beta$, drawn as
   \[
    \ba \ig{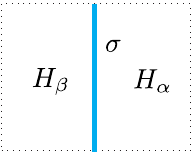} \ea \quad
    \text{or} \quad \ba \ig{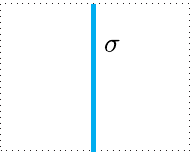} \ea \, \, \text{($H_\alpha$ and $H_\beta$
    understood)}
   \]
  The line for $\sigma$ is thick and coloured differently, so as to distinguish it from the lines for the functors $\alpha_g$ and $\beta_g$.
  \item For each $g \in G$ a unitary natural isomorphism $\sigma(g) \colon \beta_g
  \circ \sigma \Iso \sigma \circ \alpha_g$, drawn as
   \[
    \ba \ig{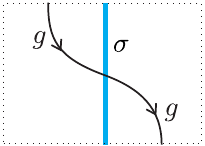} \ea \quad \ba \xymatrix{\beta_g \circ \sigma
    \ar@{=>}[d]^{\sigma(g)} \\ \sigma \circ \alpha_g} \ea
   \]
  \end{itemize}
such that
 \[
  \ba \xymatrix{ & \beta_{g_2} \circ \beta_{g_1} \circ \sigma \ar@{=>}[dl]_{\psi(g_2,g_1)
  \ast \id} \ar@{=>}[r]^{\id \ast \sigma(g_1)} & \beta_{g_2} \circ \sigma
  \circ \alpha_{g_1} \ar@{=>}[dd]^{\sigma(g_2) \ast \id} \\ \beta_{g_2g_1} \circ
  \sigma \ar@{=>}[dr]_{\sigma(g_2 g_1)} \\ & \sigma \circ \alpha_{g_2 g_1} & \sigma
  \circ \alpha_{g_2} \circ \alpha_{g_1} \ar@{=>}[l]^{\id \ast \phi(g_2,g_1)}} \ea \quad
  \ba \xymatrix{\sigma \ar@{=>}[d]_{\psi(e) \ast \id}
  \ar@{=>}[dr]^{\id \ast \phi(e)} \\ \beta_e \circ \sigma
  \ar@{=>}[r]_{\sigma(e)} & \sigma \circ \alpha_e} \ea
 \]
commute, or in string diagrams,
 \be \label{sig1}
  \ba \ig{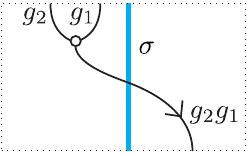} \ea = \ba \ig{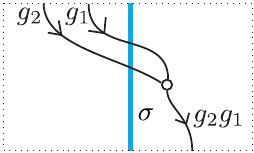} \ea \quad \text{and} \quad
  \ba \ig{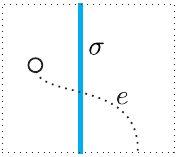} \ea = \ba \ig{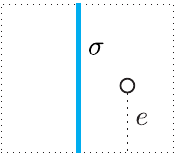} \ea.
 \ee
We will draw the inverse 2-isomorphisms $\sigma(g)^* \colon
\sigma \circ \alpha_g \Rightarrow \beta_g \circ \sigma$ as
 \[
  \ba \ig{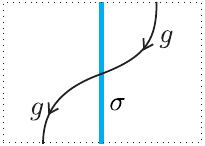} \ea \quad \ba \xymatrix{\sigma \circ \alpha_g
    \ar@{=>}[d]^{\sigma(g)^*} \\ \beta_g \circ \sigma} \ea \, .
 \]
These satisfy $\sigma(g)^* \sigma(g) = \id$ and $\sigma(g)
\sigma(g)^* = \id$, that is,
 \be \label{siginv}
 \ba \ig{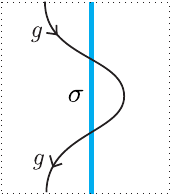} \ea = \ba \ig{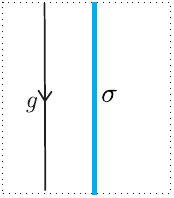} \ea \quad \text{and} \quad
 \ba \ig{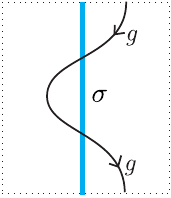} \ea = \ba \ig{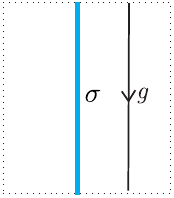} \ea \, .
 \ee
We will abbreviate all of this data $(\sigma, \{\sigma(g)\})$ simply as $\sigma$. Observe that a morphism of 2-actions of $G$, which might be called an {\em intertwiner}, really {\em does} have an `intertwining' aspect to it when expressed in terms of string diagrams.

\subsubsection{2-morphisms} Finally, if $\alpha$ and $\beta$ are unitary 2-representations of $G$, and $\sigma, \rho \colon \alpha \rightarrow \beta$ are morphisms between them, then a 2-morphism $\theta \colon \sigma \Rightarrow \rho$ is a modification from $\sigma$
 to $\rho$. Thus, $\theta$ is a natural transformation $\sigma$ to $\rho$, drawn as
 \[
  \ba \ig{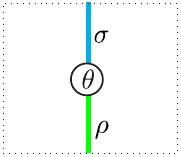} \ea \quad \ba \xymatrix{\sigma
  \ar@{=>}[d]^{\theta} \\ \rho} \ea \, ,
 \]
such that
 \[
 \xymatrix{ & \beta_g \circ \sigma \ar@{=>}[dl]_{\sigma(g)}
 \ar@{=>}[dr]^{\id \ast \theta} \\ \sigma \circ \alpha_g
 \ar@{=>}[dr]_{\theta \ast \id} & & \beta_g \circ \rho
 \ar@{=>}[dl]^{\rho(g)} \\ & \rho \circ \alpha_g}
 \]
commutes, or in string diagrams,
 \be \label{mod}
  \ba \ig{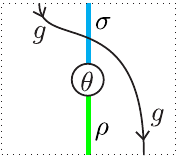} \ea = \ba \ig{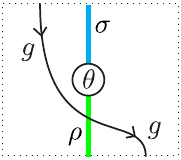} \ea \, .
 \ee
We trust that the simplicity of these diagrams has persuaded the reader the string diagrams are a useful notation for working with 2-representations. We will develop this notation further as we go along.

\subsection{Examples\label{exsec}}  We now give some examples to illustrate these ideas; we will say more about them in Section \ref{gexam} once we have established the geometric interpretation of unitary 2-representations in terms of equivariant gerbes. We urge the reader to consult \cite{ref:ganter_kapranov_rep_char_theory} for additional examples of group actions on linear categories.

\subsubsection*{2-representations can be strictified}
Before we give the examples, let us first clear up some potential confusion. A 2-representation $\alpha$ is called {\em strict} if all the
coherence isomorphisms are identities.
\begin{lem} Every 2-representation is equivalent inside $\TRep(G)$ to a strict 2-representation.
\end{lem}
\begin{proof} The proof is essentially an application of the 2-Yoneda lemma, which can be found for instance in \cite[page 60]{ref:fantechi}. Given a
2-representation of $G$ on a 2-Hilbert space $H$, we can define a corresponding strict 2-representation of $G$ on the 2-Hilbert space
 \[
  \Hom_{\TRep(G)} (\Hilb[G], \alpha),
 \]
where $\Hilb[G]$ is the category of $G$-graded Hilbert spaces on which $G$ acts by left multiplication.
\end{proof}
We warn the reader that this does not mean that there is no information in the coherence isomorphisms --- it just means that this information can always be shifted into the structure of a new and bigger category, if one wishes to do so. In other words, {\em a strict 2-representation on a `big' 2-Hilbert space} (such as `the category of all such and such') {\em is not necessarily trivial} --- to decide this, one has to calculate the 2-cocycle of the corresponding equivariant gerbe, as we explain in Section \ref{extract}.

\subsubsection{Automorphisms of groups} Suppose $G \subseteq \Aut(K)$ is a subgroup of the automorphism group of a
finite group $K$. This gives rise to a unitary 2-representation of $G$ on the 2-Hilbert space $\Rep(K)$ by precomposition. That is, if
$V$ is a unitary representation of $K$, then $g \cdot V \equiv V^g$ has the same underlying vector space
except that the action of $k \in K$ on $V^g$ corresponds to the action of $g^\mi \cdot k$ on $V$. This is of course a {\em strict} 2-representation, but it is not necessarily frivolous, as we shall see in the next example. Also note that any 2-representation of this form will necessarily be {\em unitary}, because it can only permute irreducible representations of the same dimension amongst each other.

\subsubsection{The metaplectic representation} A good example of a nontrivial 2-representation of the above sort is the action of $SL_2(\mathbb{R})$ on $\Rep(\text{Heis})$, the category of representations of the Heisenberg group. Of course, $SL_2(\mathbb{R})$ is not a {\em finite} group, but all the definitions above still apply.

The Heisenberg group arises in quantum mechanics (see for example \cite{ref:carter_segal_macdonald}). It is the 3-dimensional Lie group with underlying manifold
$\mathbb{R}^2 \times U(1)$ --- with $\mathbb{R}^2$ thought of as phase space with elements being pairs $v = (z, p)$
--- and multiplication defined by
 \[
  (v, e^{i \theta}) \cdot (w, e^{i \phi}) = (v + w, e^{i\omega (v, w)} e^{i(\theta +
  \phi)}),
 \]
where $\omega(v, w) = \frac{1}{2}(v_z w_p - v_p w_z)$ is the canonical symplectic form on $\mathbb{R}^2$. Up to
isomorphism, there is only one irreducible representation of $\text{Heis}$ on a separable Hilbert space, with
the $U(1)$ factor acting centrally. Namely, the action on $L^2(\mathbb{R})$ given by
 \[
  (z \cdot f) (x) = e^{izx} f(x), \quad (p\cdot f)(x) = f(x-p).
 \]
Since there is only one irreducible representation, $\Rep(\text{Heis})$ is a one-dimensional 2-Hilbert space.

Now $SL_2 (\mathbb{R})$ is the group of symplectomorphisms of $\mathbb{R}^2$, hence it acts as automorphisms of
Heis, giving rise to a unitary 2-representation of $SL_2(\mathbb{R})$ on $\Rep(\text{Heis})$ via the standard
prescription $(g \cdot \rho)(v, e^{i \theta}) = \rho(g^\mi \cdot v, e^{i \theta})$.

This gives rise to a nontrivial projective representation of $SL_2 (\mathbb{R})$ on $L^2(\mathbb{R})$; the fact that the projective factor cannot be removed is known as the `metaplectic anomaly'. Indeed, the viewpoint of 2-representations elucidates somewhat the nature of this anomaly.  It might seem strange at
first that the action of $SL_2(\mathbb{R})$ --- the symmetry group of the classical phase space $\mathbb{R}^2$
--- does not survive quantization, becoming instead a projective representation. However $SL_2(\mathbb{R})$
{\em does} act on $\Rep(\text{Heis})$, the collection of {\em all} quantizations. From this we see that the
`anomaly' arose from an attempt to decategorify this action, by artificially choosing a fixed quantization
$\rho$.

\subsubsection{2-representations from exact sequences\label{exactseqsec}} We've seen how an action of $G$ on another group $K$ gives rise
to a unitary 2-representation of $G$ on $\Rep(K)$. The same can be said for a `weak' action of $G$ on $K$. Suppose
 \[
  1  \rightarrow K \stackrel{i}{\hookrightarrow} E \stackrel{\pi}{\twoheadrightarrow} G  \rightarrow 1
 \]
is an exact sequence of finite groups, which has been equipped with a set-theoretic section $s \colon G \rightarrow
E$ such that $s(e) = e$. We can think of this data as a homomorphism of 2-groups
 \[
  G \rightarrow AUT(K)
 \]
where $AUT(K)$ is the 2-group whose objects are the automorphisms of $K$ and whose morphisms are given by
conjugation (see \cite{ref:baez_lauda_2-groups, ref:baez_higher_schreier_theory}). Explicitly, one thinks of the group $K$ as being the morphisms of a one-object category (also
denoted $K$), and for each $g \in G$, $g \colon K \rightarrow K$ is the functor defined by conjugating in $E$,
 \[
  g \cdot k := s(g) k s(g)^\mi
 \]
where we have identified $K$ with its image in $E$. This determines a $K$-valued 2-cocycle $\varphi$
having the property that
 \[
 g_2 \cdot g_1 \cdot k = \varphi(g_2, g_1) [ (g_2 g_1)\cdot k ] \varphi(g_2, g_1)^\mi
 \]
for all $k \in K$.

This data gives rise to a unitary 2-representation $\alpha$ of $G$ on $\Rep(K)$, by precomposition. Explicitly,
if $\rho$ is a representation of $K$, and $g \in G$, then $\alpha_g(\rho)$ has the same underlying vector space
as $\rho$, with the action of $K$ given by
 \[
  \alpha_g (\rho) (k) = \rho(g^\mi \cdot k).
 \]
The coherence natural isomorphisms $\phi(g_2, g_1) \colon \alpha_{g_2} \circ \alpha_{g_1} \Rightarrow \alpha_{g_2
g_1}$ have components
 \[
  \phi(g_2, g_1)_\rho = \rho(\varphi(g_1^\mi, g_2^\mi))
 \]
while $\phi(e) \colon \id \Rightarrow \alpha_e$ is just the identity.

\subsubsection{Other examples of 2-representations} One expects to find similar examples of unitary 2-representations of groups arising from
automorphisms of other geometric or algebraic structures --- for instance, the automorphisms of a {\em rational vertex operator algebra} or of an {\em affine lie algebra} will act on their category of representations, which in good cases are 2-Hilbert spaces.

\subsubsection{Morphisms of 2-representations from morphisms of exact sequences} We have seen how one obtains a
2-representation of $G$ from an exact sequence of groups (equipped with a set-theoretic section) with $G$ as the
final term, or equivalently from a weak action of $G$ on another group. A {\em morphism} of such a structure
gives rise to a morphism of 2-representations by induction. Indeed, suppose we have a map of exact sequences
 \[
  \ba \xymatrix{1 \ar[r] & K \ar[d]^{f_0} \ar[r] & E \ar[d]^{f_1} \ar[r] & G \ar[r] \ar[d]^\id & 1 \\
   1 \ar[r] & L \ar[r] & F \ar[r] & G \ar[r] & 1} \ea.
\]
In higher category language, this is essentially the same thing as a morphism inside the 2-category
 \[
  \Hom (G, \mathcal{G}\text{roups})
 \]
of weak 2-functors, transformations and modifications from $G$ (thought of as a one object 2-category with only
identity 2-cells) to the 2-category of groups (objects are groups, morphisms are functors, 2-morphisms are
natural transformations). Then by inducing along $f_0$ we get a map
 \[
  \sigma \equiv \Ind(f_0) \colon \Rep(K) \rightarrow \Rep(L)
  \]
and also natural isomorphisms
 \[
  \sigma(g) \colon \beta_g \circ \sigma \Rightarrow \sigma \circ \alpha_g,
 \]
where $\alpha$ and $\beta$ are the associated 2-representations of $G$ on $\Rep(K)$ and $\Rep(L)$ respectively.
In other words, {\em a map of exact sequences gives rise to a morphism of 2-representations}.

\subsection{More graphical elements\label{slem}} The reason we have been drawing arrows on the strings representing the functors $\alpha_g$ involved in a 2-representation $\alpha$ is to conveniently distinguish group elements from their inverses: if a downward pointing section of a string is labeled `$g$' then it represents $\alpha_g$, and upward pointing sections of the same string represent $\alpha_{g^\mi}$. Using this convention we now construct some new graphical elements from the old ones. From now on we drop the bounding boxes on the diagrams.

Define $\eta_g \colon \id \Rightarrow \alpha_{g^\mi} \circ \alpha_g$ and $\epsilon_g \colon \alpha_g \circ \alpha_{g^\mi} \Rightarrow \id$ as:
 \[
 \eta_g = \ba \ig{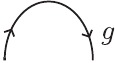} \ea := \ba \ig{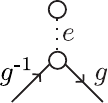} \ea  \ \ \ \equiv  \ba \xymatrix@R=0.5cm{ \id \ar@{=>}[d]^{\, \phi(e)} \\ \alpha_e
 \ar@{=>}[d]^{\, \phi(g^\mi, g)^\mi} \\ \alpha_{g^\mi} \circ \alpha_g} \ea
 \]
 \[
  \epsilon_g = \ba \ig{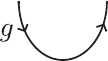} \ea := \ba \ig{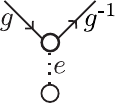}
 \ea  \ \ \ \equiv \ba \xymatrix@R=0.5cm{\alpha_g \circ \alpha_{g^\mi}
 \ar@{=>}[d]^{\phi(g, g^\mi)} \\ \alpha_e \ar@{=>}[d]^{\phi(e)^\mi} \\
 \id} \ea
 \]
These are indeed unitary natural transformations, since their inverses are clearly given by
 \[
  \eta_g^* = \ba \ig{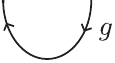} \ea := \ba \ig{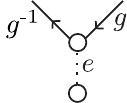} \ea \qquad
  \qquad \epsilon_g^* = \ba \ig{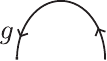} \ea := \ba \ig{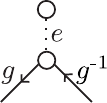}
  \ea \, .
 \]
In other words, we have the ``no
loops'' and ``merging'' rules
 \[
 \ba \ig{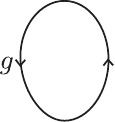} \ea = \ba \ig{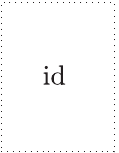} \ea \qquad \qquad \ba
 \ig{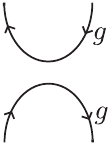} \ea = \ba \ig{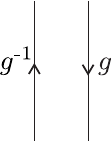} \ea
 \]
and similarly for the reverse orientations.

We now show that these new graphical elements behave as their string diagrams suggest. The first part of the following
lemma actually says in more orthodox terminology that `$\alpha_g$ is an ambidextrous adjoint equivalence from the underlying 2-Hilbert space to itself', or more precisely `for all $g \in G$, $\alpha_g \dashv \alpha_{g^\mi}$ via $(\eta_g,
\epsilon_g)$'. But it's the simple fact that these string diagrams can be manipulated in the obvious intuitive
fashion which is more important for us here.
\begin{lem} \label{ambilem} Suppose $\alpha$ is a 2-representation of $G$. The following graphical moves hold:
 \[
 \text{ \em (i)}
 \ba \ig{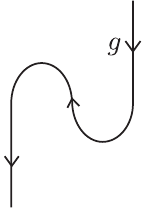} \ea = \ba \ig{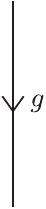} \ea = \ba \ig{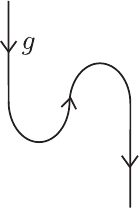} \ea
 \qquad \text{\em (ii)} \quad \ba \ig{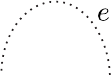} \ea = \ba \ig{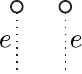} \ea
 \]
 \[
  \text{\em (iii)} \ba \ig{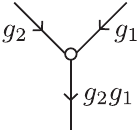} \ea = \ba \ig{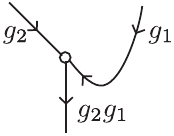} \ea = \ba \ig{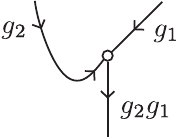} \ea  \qquad
  \text{\em (iv)} \ba \ig{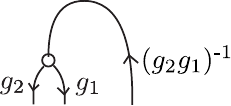} \ea = \ba \ig{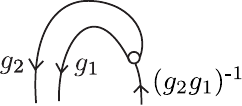} \ea .
 \]
\end{lem}
\begin{proof} (i) The first equation as proved as follows,
 \[
 \ba \ig{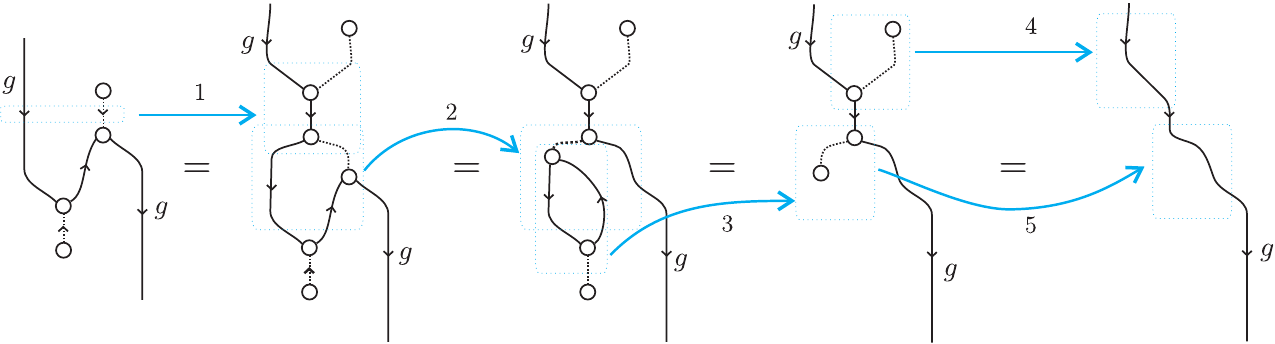} \ea .
 \]
In step 1 we zip together using the rule (\ref{zr1}a), in 2 we slide the button around using (\ref{br1}b), in 3
we unzip again using (\ref{zr1}a) and in 4 and 5 we contract the identity string using (\ref{br2}b). The other
equations are proved similarly.
 \end{proof}
Now we record for further use some allowable graphical manipulations for {\em morphisms} of 2-representations.
 \begin{lem} \label{onestlem} Suppose $\sigma \colon \alpha \rightarrow \beta$ is a morphism of 2-representations. The following graphical
 moves hold:
 \[
  \text{(i)} \ba \ig{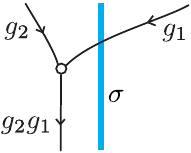} \ea = \ba \ig{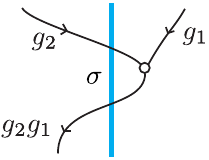} \ea \qquad \text{(ii)} \ba \ig{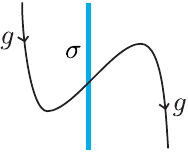} \ea = \ba
  \ig{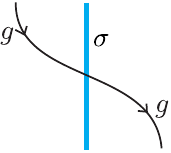} \ea
  \]
 \end{lem}
 \begin{proof} (i) is proved as follows,
  \[
  \ba \ig{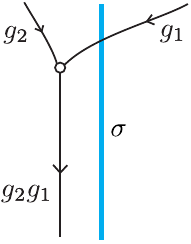} \ea \stackrel{\text{(a)}}{=} \ba \ig{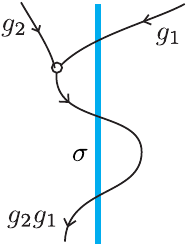} \ea
  \stackrel{\text{(b)}}{=} \ba \ig{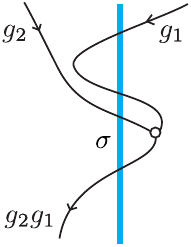} \ea
  \stackrel{\text{(c)}}{=} \ba \ig{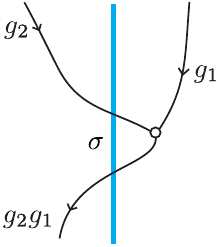} \ea \, ,
  \]
 where (a) uses the inverse rule (\ref{siginv}a), (b) uses the button-dragging rule
 (\ref{sig1}a), and (c) uses the reverse inverse rule
 (\ref{siginv}b).
 \end{proof}

\subsection{The even-handed structure on 2$\mathcal{H}\text{ilb}$\label{EvenSec}}
In order to make the 2-character functorial, we will need to have tight control on the ambidextrous (simultaneous left and right) adjoints in $\THilb$. This is accomplished via an {\em even-handed structure} --- a coherent system for turning right adjoints into left adjoints. We first define this notion for general 2-categories, and then we show that $\THilb$ has a canonical such structure.

\subsubsection{Even-handed structures on general 2-categories}
Let us first be clear about our terminology.
\begin{defn}
 An {\em ambidextrous adjoint} of a morphism $F \colon A \rightarrow B$ in a 2-category is a quintuple $\langle F^* \rangle \equiv (F^*, \eta, \epsilon, n, e)$ where  $F^* \colon B \rightarrow A$ is a morphism, $\eta \colon \id_A \Rightarrow F^*F$ and $ \epsilon \colon FF^* \Rightarrow \id_B$ are unit and counit maps exhibiting $F^*$ as a right adjoint of $F$, and $n \colon \id_B \Rightarrow F F^*$ and $e \colon F^* F \Rightarrow \id_A$ are unit and counit maps which exhibit $F^*$ as a left adjoint of $F$.
 \end{defn}
We write the data of a particular ambidextrous adjoint of $F$ in string diagrams as
 \[
   \langle F^* \rangle  = \left( \ba \ig{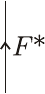} \ea, \ba \ig{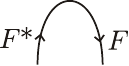} \ea, \ba
   \ig{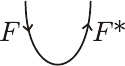} \ea, \ba \ig{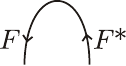} \ea, \ba \ig{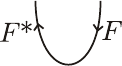} \ea \right) \, .
 \]
We can organize the choices of ambidextrous adjoints for $F$ into a groupoid $\Amb(F)$ which we call the {\em ambijunction groupoid} of $F$, as follows. An object is a choice of ambidextrous adjoint $\langle F^* \rangle$ of $F$. A morphism $\gamma \colon \langle F^* \rangle \rightarrow  \langle (F^*)' \rangle$ of ambidextrous adjoints of $F$ is an invertible 2-morphism $\gamma \colon  F^* \Rightarrow (F^*)'$, drawn as
 \[
  \ba \ig{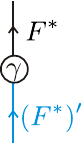} \ea \,
 \]
such that `twisting' the unit and counits of $\langle F^* \rangle$ by $\gamma$ results in $\langle (F^*)' \rangle$, that is,
 \[
  \left( \ba \ig{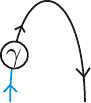} \ea, \ba \ig{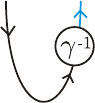} \ea, \ba
  \ig{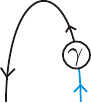} \ea, \ba \ig{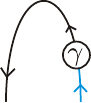} \ea \right) = \left( \ba
  \ig{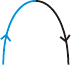} \ea, \ba \ig{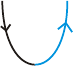} \ea, \ba \ig{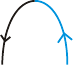} \ea, \ba
  \ig{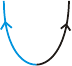} \ea \right).
 \]
We write $[\Amb(F)]$ for the set of isomorphism classes in the ambijunction groupoid of $F$, and we write the class of a particular ambidextrous adjoint $\langle F^* \rangle$ as $[F^*]$. The properties of the ambijunction groupoid are summarized in the following elementary but important lemma.
\begin{lem} \label{amblem1} Suppose $F$ is a morphism in a 2-category and that the groupoid $\Amb(F)$ is nonempty. Then:
 \begin{enumerate}
  \item There is at most one arrow between any two ambidextrous adjunctions in $\Amb(F)$.
  \item The group $\Aut(F)$ of automorphisms of $F$ acts freely and transitively on $[\Amb(F)]$ by twisting the unit and counit maps which display $F^*$ as a left adjoint of $F$,
      \[
  [F^*] \stackrel{\alpha}{\mapsto} \left[ \ba \ig{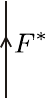} \ea, \ba \ig{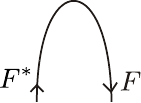} \ea, \ba \ig{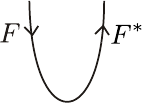} \ea, \ba
  \ig{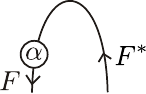} \ea, \ba \ig{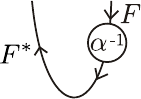} \ea \right].
 \]
 \end{enumerate}
\end{lem}
\noindent Now suppose that $\theta \colon  F \Rightarrow G$ is a 2-morphism, drawn as
\[
 \ba \ig{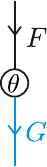} \ea,
 \]
and that choices $\langle F^* \rangle$, $\langle G^* \rangle$ of ambidextrous adjoints of $F$ and $G$ have been made. The {\em right and left daggers} of $\theta$ are defined to be the 2-morphisms $\theta^\dagger, ^\dagger\!\theta \colon  G^*
\Rightarrow F^*$ given by
 \be \label{adjoints1}
  \theta^\dagger := \ba\ig{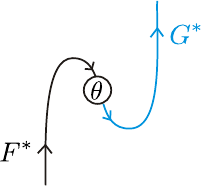} \ea \quad \quad\ ^\dagger\!\theta
  := \ba \ig{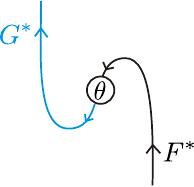} \ea.
 \ee
In other words, the right dagger $\theta^\dagger$ is constructed from the data of $F^*$ and $G^*$ being {\em right} adjoints of $F$ and $G$ respectively, while the left dagger $^\dagger\!\theta$ is constructed from the data of $F^*$ and $G^*$ being {\em left} adjoints of $F$ and $G$ respectively. It is clear that $\theta^\dagger$ need not equal $^\dagger\!\theta$ because they transform differently under the action of automorphisms of $F$ and $G$. Indeed, if we use automorphisms $\alpha_F\colon F \Rightarrow F$ and $\alpha_G\colon G \Rightarrow G$ to twist the left adjoint unit and counit maps, we see that the right daggers remain unchanged while the left daggers transform as
\[
 \ba \ig{e30.pdf} \ea \quad \mapsto \quad \ba \ig{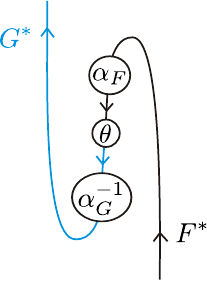} \ea.
 \]
 However, the question of whether the right dagger $\theta^\dagger$ is equal to the left dagger $^\dagger\!\theta$ only depends on the {\em isomorphism classes} of $ \langle F^* \rangle$ and $\langle G^* \rangle$ in their respective ambijunction groupoids, as the reader will verify by an elementary string diagram calculation.

This suggests the following definition. Firstly, observe that one can compose ambidextrous adjoints in the obvious way, and also that every identity morphism $\id_A$ in a 2-category has the trivial ambidextrous adjunction associated to it, with unit maps $\id_A \Rightarrow \id_A \circ \id_A$ given by the unit isomorphisms and counit maps $\id_A \circ \id_A \Rightarrow \id_A$ given by their inverses. Also, we say that a 2-category $\mathcal{C}$ {\em has ambidextrous adjoints} if every morphism has an ambidextrous adjoint. Note that having ambidextrous adjoints is a {\em property} of $\mathcal{C}$ --- we do {\em not} require that permanent choices of these adjoints have been made from the start.
\begin{defn} An {\em even-handed structure} on a 2-category with ambidextrous adjoints is a choice $F^{[*]} \in [\Amb(F)]$ of isomorphism class of ambidextrous adjoint for every morphism $F$, such that:
 \begin{enumerate}
  \item $\id^{[*]}$ is the class of the trivial ambidextrous adjunction for every identity morphism,
  \item $(G \circ F)^{[*]} = F^{[*]} \circ G^{[{*}]}$ for all composable pairs of morphisms, and
  \item $\theta^\dagger = {} ^\dagger\!\theta$ for every 2-morphism $\theta \colon  F \Rightarrow G$, provided they are computed using ambidextrous adjoints from the classes $F^{[*]}$ and $G^{[*]}$.
  \end{enumerate}
An {\em even-handed 2-category} is a 2-category with ambidextrous adjoints equipped with an even-handed structure. \end{defn}
In other words, instead of stipulating a specific ambidextrous adjoint $\langle F^* \rangle$ for every morphism $F$, an even-handed structure selects only an {\em isomorphism class} $F^{[*]}$ of such adjoints, in such a way that the resultant choices are compatible with composition and ensure that the left and right daggers of 2-morphisms always agree. The idea is that calculations involving ambidextrous adjoints usually do not depend on the actual ambidextrous adjoints themselves, but they will depend on their isomorphism classes, and so these classes need to be given as extra information.  We will say that a particular choice of ambidextrous adjoint $\langle F^* \rangle$ is {\em even-handed} if it is a member of the specified class $F^{[*]}$.

\subsubsection{Geometric interpretation of even-handed structures}
\begin{figure}[t]
 \centering
 \ig{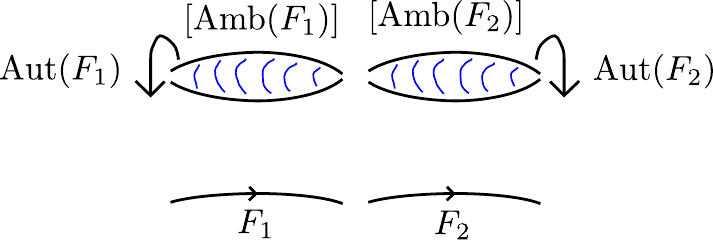}
 \caption{\label{ambgerbe} The `ambijunction gerbe' of a 2-category with ambidextrous adjoints.}
 \end{figure}

We encourage the reader to think of an even-handed structure geometrically in the following way. If $\mathcal{C}$ is a 2-category with ambidextrous adjoints, the sets $[\Amb(F)]$ of isomorphism classes of ambidextrous adjoints for each morphism $F$ can be thought of as forming a `gerbe-like' structure
\[
 \Amb(\mathcal{C}) \rightarrow \text{1-Mor}(\mathcal{C})
\]
over the morphisms in $\mathcal{C}$, which we call the {\em ambijunction gerbe} (see Figure \ref{ambgerbe}). An even-handed structure is then an `even-handed trivialization' of this gerbe. This is analogous to Murray and Singer's reformulation of the notion of a spin structure on a Riemannian manifold $M$ as a trivialization of the {\em spin gerbe} \cite{ref:murray_singer}. We develop this geometric analogy further in \cite{ref:bartlett_even_handed}; for now we wish to stress that a given 2-category can admit many different even-handed structures.

\subsubsection{Diagram manipulations requiring an even-handed structure}
The following lemma illustrates the sort of diagram manipulations which require an even-handed structure; we will need these results later. We use the convention that in the context of an even-handed 2-category, all string diagram equations involving computations with ambidextrous adjoints are to be interpreted as being appended with the disclaimer `provided the ambidextrous adjoints are chosen from the classes stipulated by the even-handed structure'; in other words they hold if the chosen ambidextrous adjoints are even-handed.
\begin{lem} \label{evenhandlem} For a 2-morphism $\theta \colon  F \Rightarrow G$ in an even-handed 2-category, the
following equations, together with all obvious variations, hold:
 \[
  \text{(i)} \, \, \ba \ig{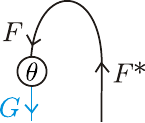} \ea = \ba \ig{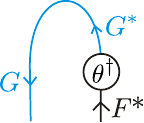} \ea \quad \quad \text{(ii)} \, \,\ba \ig{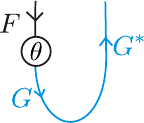} \ea = \ba
  \ig{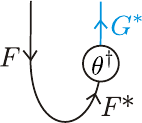} \ea
 \]
 \[
 \text{(iii) } \, \, (\theta^\dagger)^\mi = (\theta^\mi)^\dagger \text{ (when $\theta$ is
  invertible).}
  \]
 \end{lem}
\begin{proof} For {\em both} (i) and (ii) to hold, we need the ambidextrous adjoints to be even-handed. For instance, (i) is proved as follows:
\[
  \ba \ig{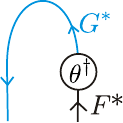} \ea = \ba \ig{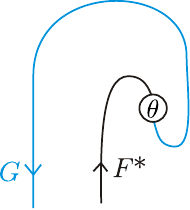} \ea = \ba \ig{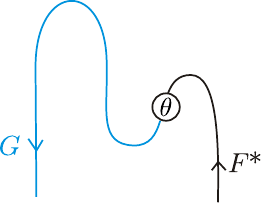} \ea = \ba \ig{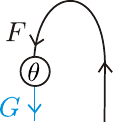} \ea.
\]
The rest are proved similarly.
\end{proof}

\subsubsection{The even-handed structure on 2$\mathcal{H}\text{ilb}$}
To discuss adjoints and even-handed structures in the concrete setting where the 2-category consists of categories, functors and natural transformations, it is convenient to change variables from the unit and counit natural transformations to the associated adjunction isomorphisms
 \[
  \varphi \colon  \Hom(Fx, y) \stackrel{\cong}{\rightarrow} \Hom(x, F^*y).
 \]
Recall the translation between these two pictures --- if $f \colon  Fx \rightarrow y$ then $\varphi(f) = F^*(f) \circ
\eta_x$. Suppose we write
 \[
 \Adj(F \dashv F^*) = \{ \varphi \colon \Hom(Fx, y) \stackrel{\cong}{\rightarrow} \Hom(x, F^*y), \text{ natural in $x$ and $y$}\}
 \]
for the set of ways in which an adjoint pair of functors can be expressed as adjoints of each other (recall that this set is a torsor for the group $\Aut(F)$). In this language, an even-handed structure becomes a collection of bijective maps
\[
 \Psi_{F, F^*} \colon \Adj(F \dashv F^*) \rightarrow \Adj(F^* \dashv F),
 \]
one for each pair of functors capable of being expressed as adjoints of each other, which is compatible with natural isomorphisms and composition and such that $\theta^\dagger = {} ^\dagger \!\theta$ for all natural transformations $\theta$.

As an example, consider the 2-category $\THilb$ of 2-Hilbert spaces. Since 2-Hilbert spaces are semisimple, every 1-morphism $F \colon H \rightarrow H'$ in $\THilb$ has an ambidextrous adjoint. That is because each functor $F$ is freely determined up to isomorphism by the matrix of nonnegative integers $\dim \Hom(e_\mu, Fe_i)$, where $e_\mu$ and $e_i$ run over a choice of simple objects in $H$ and $H'$ respectively, and so one can choose a functor $F^* \colon H' \rightarrow H$ whose associated matrix is the transpose of that of $F$. Equipping $F^*$ with the structure of an ambidextrous adjoint $\langle F^* \rangle$ amounts to freely (and apriori {\em independently}) making a choice of linear isomorphisms
 \begin{align*}
  \phi_{i, \mu} \colon \Hom(Fe_i, e_\mu) & \rightarrow \Hom(e_i, F^*e_\mu) \\
   \intertext{and}
   \psi_{\mu, i} \colon \Hom(F^*e_\mu, e_i) & \rightarrow \Hom(e_\mu, F e_i)
 \end{align*}
respectively. An even-handed structure is precisely a system which pairs these isomorphisms $\phi$ and $\psi$ together, so that knowing one determines the other. To repeat: $\THilb$ does not have canonically given ambidextrous adjoints. Rather, what {\em is} canonical (as we state in the following proposition) is the function which {\em turns} right adjoints into left adjoints. Readers familiar with these ideas in an algebraic geometry context might want to translate this into the statement that {\em every 2-Hilbert space comes canonically equipped with a trivial Serre functor}. The following is shown in \cite{ref:bartlett_even_handed}.
\begin{prop} The 2-category $\THilb$ comes equipped with a canonical even-handed structure, given at the level of adjunction isomorphisms by sending
 \[
  \varphi \mapsto * \, \varphi^* \, *
 \]
where $\varphi^*$ is the adjoint of $\varphi$ in the ordinary sense of maps between Hilbert spaces.
\end{prop}
\noindent That is, the even-handed structure sends
  \[
  \ba \xymatrix @C=0.25in {\Hom(Fx, y) \ar[r]^-\varphi & \Hom(x, F^*y)} \ea \quad \mapsto \quad \ba \xymatrix{\Hom(F^*y, x) \ar[d]_{*} & \Hom(y, Fx)  \\
  \Hom(x,F^*y) \ar[r]_{\varphi^*} & \Hom(Fx, y)\ar[u]_{*} }  \ea.
 \]
This formula resembles the formula for the adjoint of the derivative operator on a Riemannian manifold, $d \mapsto * \, d \, *$. Also note that we needed the inner products and the $*$-structure on the hom-sets for this to work --- there is no canonical even-handed structure on Kapranov and Voevodsky's 2-category $\TVect$ of 2-vector spaces, for example. Indeed, we show in \cite{ref:bartlett_even_handed} that an even-handed structure on $\TVect$ is the same thing, up to a global scale factor, as an assignment of a nonzero {\em scale factor} to each simple object --- which is precisely the extra information available in a  2-Hilbert space. This underscores the fact that an even-handed structure on a 2-category can be thought of as equipping each object of the 2-category with a `metric'.

\subsection{Even-handedness and unitary 2-representations} Consider the category $\Rep(G)$ of unitary representations of a group $G$. If $\sigma \colon \rho_1 \rightarrow \rho_2$ is an intertwining map between two unitary representations, then the adjoint $\sigma^* \colon \rho_2 \rightarrow \rho_1$ is also an intertwining map, because
 \begin{align*}
  \sigma^* \circ \rho_2 (g) &= \sigma^* \circ \rho_2^*(g^\mi) \\
   &= (\rho_2 (g^\mi ) \circ \sigma)^* \\
   &= (\sigma \circ \rho_1 (g^\mi ))^* \\
   &= \rho_1(g) \circ \sigma^*.
  \end{align*}
We will need the corresponding result for unitary 2-representations. Suppose that $\sigma \colon \alpha \rightarrow \beta$ is a morphism of unitary 2-representations of $G$. If
 \[
  \sigma^* \colon H_\beta \rightarrow H_\alpha
 \]
is a linear $*$-functor which is adjoint to the underlying functor $\sigma \colon H_\alpha \rightarrow H_\beta$,
can we equip $\sigma^*$ with the structure of a morphism of 2-representations
 \[
 \sigma^* \colon \beta \rightarrow \alpha
 \]
in such a way that $\sigma^*$ is adjoint to $\sigma$ in $\TRep(G)$? The answer is {\em yes} --- because $\alpha$ and $\beta$ are unitary 2-representations, and this means they are compatible with the even-handed structure on $\THilb$. A 2-representation $\alpha$ is called {\em even-handed} if for each $g \in G$
\[
  \alpha_g^{[*]} =
  \left[ \ba \ig{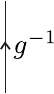} \ea, \ba \ig{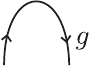} \ea,
  \ba \ig{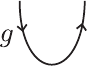} \ea, \ba \ig{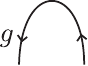} \ea, \ba \ig{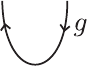} \ea \right],
\]
and one can check that unitary 2-representations indeed have this property. This means we can choose $(\sigma \circ \alpha_{g^\mi})^*$ to be $\alpha_g \circ \sigma^*$, and similarly for $\beta$. Thus we can define $\sigma^* (g) \colon \alpha_g \circ \sigma^* \Rightarrow \sigma^* \circ \beta_g$ to be
 \[
  \alpha_g \circ \sigma^* = (\sigma \circ \alpha_{g^\mi})^* \stackrel{\sigma(g^\mi)^\dagger}{\Rightarrow}
  (\beta_{g^\mi} \circ \sigma)^* = \sigma^* \circ \beta_g.
 \]
In string diagrams,
 \[
 \sigma^*(g) \stackrel{\text{drawn as}}{\equiv} \ba \ig{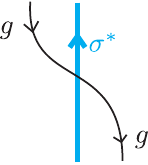} \ea := \ba \ig{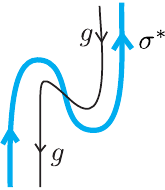} \ea = \ba \ig{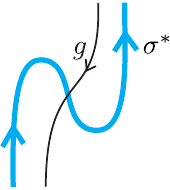} \ea
 \]
where the last simplification step uses Lemma \ref{onestlem} (ii). Since we are using even-handed adjunctions,
it would not have made a difference if we had opted to use the left dagger $^\dagger\!\sigma(g^\mi)$ instead. Written in string diagrams, this means
 \[
 \sigma^*(g) = \ba \ig{m117.pdf} \ea = \ba \ig{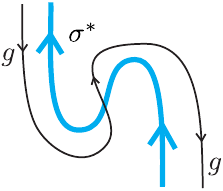} \ea.
 \]
The fact that the maps $\sigma^*(g)$ satisfy the coherence equations for a morphism of 2-representations follows routinely from the fact that $\sigma(g)$ satisfies them. In this way we have lifted the morphism $\sigma^* \colon H_\beta \rightarrow H_\alpha$ in
$\THilb$ to a morphism $\sigma^* \colon \beta \rightarrow \alpha$ in $\TRep(G)$.

We will need the following graphical moves, which formally may be regarded as establishing that the unit and
counit natural transformations for the adjunctions $\sigma \dashv \sigma^*$ and $\sigma^* \dashv \sigma$ satisfy the coherence equation for a 2-morphism in $\TRep(G)$, so that we have indeed succeeded in lifting the ambidextrous adjunction $\sigma \dashv \sigma^* \dashv \sigma$ in $\THilb$ to an ambidextrous adjunction $\sigma \dashv \sigma^* \dashv \sigma$ in $\TRep(G)$.

\begin{lem} \label{invproof} With this definition of $\sigma^*(g)$, the following equations hold:
  \begin{align*}
  &\text{(i)}\;\; \ba \ig{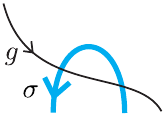} \ea = \ba \ig{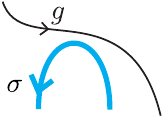} \ea &  &\text{(ii)} \;\; \ba \ig{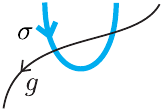} \ea = \ba \ig{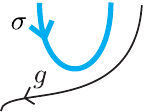} \ea \\
   &\text{(iii)} \;\; \ba \ig{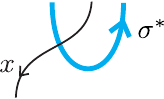}\ea = \ba \ig{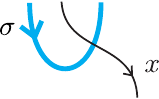} \ea &  &\text{(iv)} \;\; \sigma^*(g)^\mi = \ba \ig{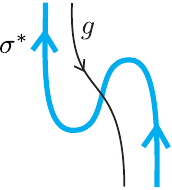} \ea
  \end{align*}
\end{lem}

\section{2-characters of 2-representations\label{charsecc}}
In this section we define the 2-character of a 2-representation. This notion was defined independently by Ganter and Kapranov \cite{ref:ganter_kapranov_rep_char_theory} while we were working on this paper. What is new in our treatment is that we show how 2-characters look especially simple when expressed in terms of string diagrams, but more importantly we use the canonical even-handed structure on $\THilb$ to show how the 2-character can be made {\em functorial} with respect to morphisms of 2-representations, as we explained in the introduction. This prepares the way for us to show that the 2-character corresponds to the `geometric character' of the associated equivariant gerbe, and also for us to show that the complexified 2-character is unitarily fully faithful.

\subsection{Definition and basic properties \label{defsec}}
\subsubsection*{2-traces}
The basic idea of 2-characters, as we explained in the introduction, is that they categorify the notion of the character of an ordinary representation of a group. Ordinary characters are defined by taking traces, so we first need to define {\em 2-traces} (Ganter and Kapranov called this the {\em categorical trace}).

\begin{defn} The {\em 2-trace} of a linear endofunctor $F \colon H \rightarrow H$ on a 2-Hilbert space $H$ is the Hilbert space
 \[
 \ttr{F} = \Nat(\id_H, F) = \left\{ \ba \ig{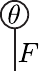} \ea \right\}.
 \]
 \end{defn}
If one thinks of $F$ via its associated matrix of Hilbert spaces $\Hom(e_j, Fe_i)$ where $e_i$ runs over a choice of simple objects for $H$, then the 2-trace corresponds to the direct sum of the Hilbert spaces along the diagonal, because a natural transformation $\id \Rightarrow F$ is freely and uniquely determined by its behaviour on the simple objects. Also recall that the inner product on $\ttr{F}$ is defined as
 \[
   \langle \theta, \theta'\rangle = \sum_{i} k_i (\theta_{e_i}, \theta'_{e_i})
 \]
where $e_i$ runs over a choice of simple objects for $H$, and $k_i = (\id_{e_i}, \id_{e_i})$ as always.

\subsubsection*{The loop groupoid}
\begin{figure}[t]
\centering
\ig{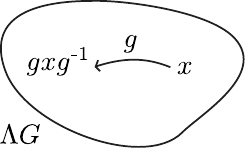}
\caption{\label{lgfig}The loop groupoid of a finite group}.
\end{figure}

In general the {\em loop groupoid} $\Lambda \mathcal{G}$ of a finite groupoid $\mathcal{G}$ is the category of functors and natural transformations from the group of integers to $\mathcal{G}$ (see \cite{ref:simon}). A special case is the loop groupoid $\Lambda G$ of a finite group $G$ which depicts the action of the group on itself by conjugation, since the objects can be identified with the elements $x \in G$, and the morphisms can be written as $gxg^\mi \sla{g} x$ (see Figure \ref{lgfig}).

\subsubsection*{The definition}The loop groupoid is particularly convenient when discussing characters. The fact that the ordinary character $\chi_\rho$ of a representation $\rho$ of $G$ is conjugation invariant can be expressed by saying it is an invariant map
 \[
  \chi_\rho \colon \Ob \Lambda G \rightarrow \mathbb{C}.
 \]
Similarly the 2-character $\chi_\alpha$ of a unitary 2-representation $\alpha$ will produce a {\em unitary equivariant vector bundle} over the group, that is, a unitary representation of the loop groupoid,
 \[
  \chi_\alpha \colon \Lambda G \rightarrow \Hilb.
 \]
We write the category of unitary equivariant vector bundles over $G$ as $\Hilb_G(G)$.

\begin{defn} The {\em 2-character} $\chi_\alpha$ of a unitary 2-representation $\alpha$ of $G$ is the unitary equivariant vector bundle over $G$ given by
 \[
  \begin{split}
   \chi_F(x)  = \ttr{\alpha_x} &= \left\{ \ba\ig{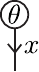} \ea \right\} \\
    \chi_\alpha(gxg^\mi \sla{g} x)\left( \ba \ig{c3.pdf}\ea \right) &=
    \ba\ig{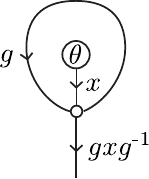} \ea.
  \end{split}
  \]
  \end{defn}
Let us verify that this definition makes sense.
\begin{prop} The 2-character $\chi_\alpha$ is indeed a unitary equivariant vector bundle over the group.
\end{prop}
\begin{proof} Using our graphical rules from Lemma \ref{ambilem}, we have
\begin{align*}
\chi_\alpha ( \sla{g_2} g_1 x g_1^\mi) \chi_\alpha (\sla{g_1} x) \left ( \ba
\ig{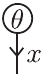} \ea \right) &= \ba \ig{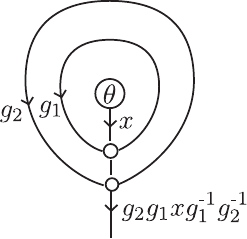} \ea =  \ba \ig{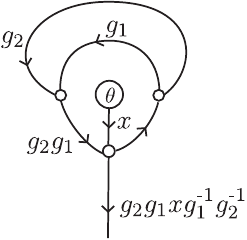} \ea \\ &= \ba\ig{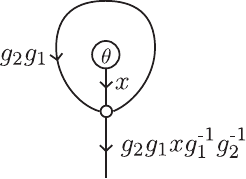} \ea
= \chi_\alpha( \la{g_2g_1} x)\left( \ba \ig{y25.pdf} \ea \right)
\end{align*}
and also
\[
 \chi_\alpha ( \sla{e} x) \left( \ba \ig{y25.pdf} \ea \right) = \ba \ig{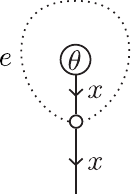} \ea = \ba \ig{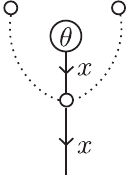} \ea =
 \ba \ig{y25.pdf} \ea.
\]
This is a {\em unitary} vector bundle because all the maps involved in its definition are unitary (see the proof of Lemma \ref{transglem} for an explicit formula).
\end{proof}

\subsection{Functoriality of the 2-character \label{funcchar}}
In this subsection we combine all the technology we have developed so far to define how to take the 2-character of a {\em morphism} of unitary 2-representations so as to obtain a morphism of the corresponding equivariant vector bundles over $G$.
\begin{defn} \label{2charmorphisms} If $\sigma \colon \alpha \rightarrow \beta$ is a morphism of unitary 2-representations, we define $\chi(\sigma) \colon \chi_\alpha \rightarrow \chi_\beta$ as the map of equivariant vector bundles over $G$ whose component at $x \in G$ is given by
 \be \label{stdiagmor}
 \begin{aligned}
  \chi(\sigma)_x \colon \chi_\alpha(x) & \rightarrow \chi_\beta (x) \\
  \ba \ig{y25.pdf} \ea & \mapsto \ba \ig{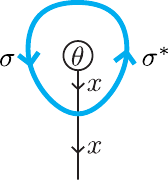} \ea
  \end{aligned}
 \ee
where $\langle \sigma^* \rangle$ is any even-handed ambidextrous adjoint for $\sigma$.
 \end{defn}
We now verify that this definition makes sense. Firstly note that it doesn't matter how we resolve the right hand
side since by Lemma \ref{invproof} (iii) we have
 \[
  \ba \ig{y29.pdf}\ea = \ba \ig{y30.pdf} \ea \, .
 \]
We write $[\TRep(G)]$ for the {\em homotopy category} of $\TRep(G)$ --- its objects are unitary 2-representations and its morphisms are isomorphism classes of 1-morphisms in $\TRep(G)$.
\begin{thm} In the situation above, the map $\chi(\sigma) \colon \chi_\alpha \rightarrow \chi_\beta$
 \begin{enumerate}
  \item does not depend on the choice of even-handed adjoint $\langle \sigma^* \rangle \in \sigma^{[*]}$,
  \item is indeed a morphism of equivariant vector bundles over $G$,
  \item does not depend on the isomorphism class of $\sigma$,
  \item is functorial with respect to composition of 1-morphisms in $\TRep(G)$,
  \end{enumerate}
 and hence $\chi$ descends to a functor
  \[
   \chi \colon [\TRep(G)] \rightarrow \Hilb_G (G).
  \]
 \end{thm}
 \begin{proof} (i) If $\langle (\sigma^*)' \rangle$ is another even-handed ambidextrous adjoint for $\sigma$, then by definition there is  a 2-isomorphism $\gamma \colon \sigma^* \Rightarrow (\sigma^{*})'$ having the property that the second equality below is valid:
   \[
    \ba \ig{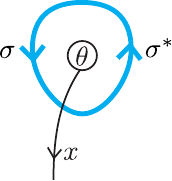} \ea = \ba \ig{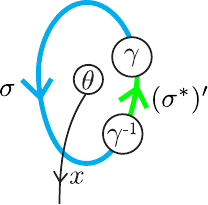} \ea \; = \ba \ig{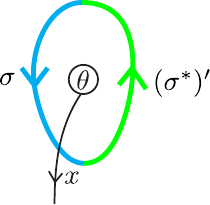} \ea.
   \]

   (ii) Using the graphical rules in Lemmas \ref{ambilem} and \ref{invproof}, we calculate:
   \[ \begin{aligned}
    \chi_\beta (\sla{g} x) \chi(\sigma)_x \left( \ba \ig{y25.pdf} \ea \right) = \ba \ig{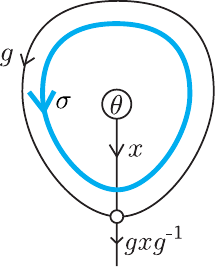} \ea = \ba
    \ig{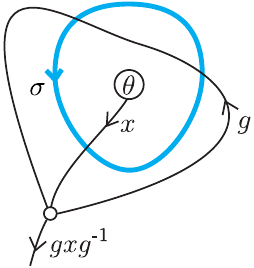} \ea = \ba \ig{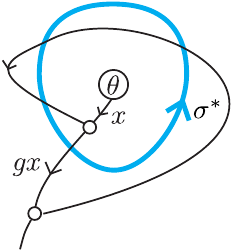} \ea \\ = \ba \ig{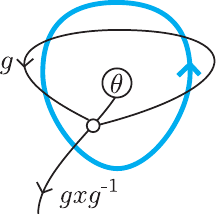} \ea = \ba \ig{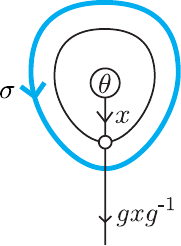} \ea = \chi(\sigma)_x
    \chi_\alpha(\sla{g} x) \left( \ba \ig{y25.pdf} \ea\right).
    \end{aligned}
   \]
   (iii) Suppose $\gamma \colon \sigma \Rightarrow \rho$ is an invertible 2-morphism in $\TRep(G)$. Then
 \[
    \ba \ig{y31.pdf} \ea = \ba \ig{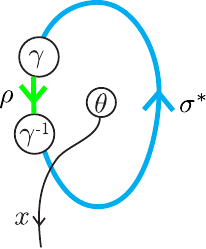} \ea = \ba \ig{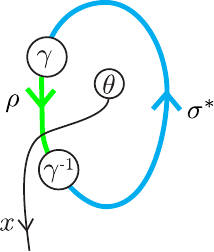} \ea = \ba \ig{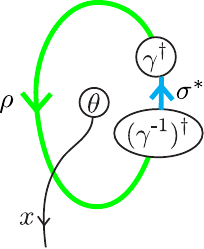} \ea = \ba \ig{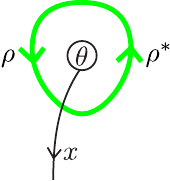}
    \ea
    \]
 where the second step uses the fact that $\gamma^\mi$ is a 2-morphism in $\TRep(G)$, while the last two steps use Lemma \ref{evenhandlem}.

  (iv)  Suppose $\delta \colon \beta \rightarrow \gamma$ is another 1-morphism in $\TRep(G)$. Since an even-handed structure respects composition, we can choose $\langle (\delta \circ \sigma)^* \rangle$ to be $\langle \sigma^* \rangle \circ \langle \delta^* \rangle$, thus $\chi(\delta \circ \sigma) = \chi(\rho) \circ \chi(\sigma)$. Similarly $\chi(\id) = \id$, because $\id^{[*]}$ is by definition the trivial ambidextrous adjunction.
   \end{proof}
We hope that these diagrammatic proofs have convinced the reader of the utility of the string diagram notation. Explicit {\em expansions} of these diagrams in terms of concrete formulas can be found in Section \ref{2char2rep}. We remark here that the behaviour of the 2-character on morphisms really does use the canonical even-handed structure on $\THilb$ in an essential way. That is, if $\sigma \colon \alpha \rightarrow \beta$ is a morphism of 2-representations, then $\chi(\sigma) \colon \chi_\alpha \rightarrow \chi_\beta$ {\em depends} on the scale factors $k_i$ and $k_\mu$ on the simple objects in $H_\alpha$ and $H_\beta$. The easiest way to see this is from our main theorem which states that the 2-character of a unitary 2-representation corresponds to the geometric character of its associated equivariant gerbe --- and the behaviour of the geometric character on morphisms of equivariant gerbes really does depend on the metrics on the gerbes (see Section \ref{fgerbes}).

\section{Equivariant gerbes\label{gggsec}}
In this section we introduce the principal geometric actors of this paper --- finite equivariant gerbes equipped with metrics, the 2-category which they constitute, and the twisted character of an equivariant vector bundle over a gerbe. As we mentioned in the introduction, we encourage the reader to think of a finite equivariant gerbe equipped with a metric as the categorification, in our finite discrete setting, of an {\em equivariant hermitian line bundle} over a compact complex manifold equipped with a metric. Our main aim in this section is to introduce the necessary geometric language so as to be able to state Theorem \ref{Simthm} below (a result due to Willerton \cite{ref:simon}) which {\em links} these two pictures --- it says that the twisted character map gives a unitary isomorphism from the space of isomorphism classes of equivariant vector bundles over an equivariant gerbe to the space of flat sections of a certain line bundle. After we have established the correspondence between 2-representations and equivariant gerbes in Section \ref{togsec}, this result will imply that the complexified 2-character functor is unitarily fully faithful.

\subsection{The definition} In this subsection we define discrete equivariant gerbes. Our definition agrees with that of Behrend and Xu \cite{ref:behrend_xu} if one specializes their notion to this simplified setting --- though we also add in the idea of a {\em metric}.

\subsubsection*{$U(1)$-torsors and their tensor products} A $U(1)$-torsor is a set with a free and transitive left
action of $U(1)$. The tensor product $P \otimes Q$ of two $U(1)$-torsors is the torsor obtained from the
cartesian product $P \times Q$ by identifying $(e^{i \theta} p, q)$ with $(p, e^{i \theta} q)$ for any $p \in
P$, $q \in Q$ and $e^{i \theta} \in U(1)$; the equivalence class of $(p,q)$ is denoted $p \otimes q$.

\subsubsection*{Equivariant gerbes \label{EqGerbeSec}} Let $X$ be a left $G$-set; we will only deal with {\em discrete} $G$-sets so we write the elements of $X$ as $i,j,k$ etc. We think of $X$ via its associated {\em action groupoid} $X_G$, which has objects the elements $i \in X$ and morphisms $g \cdot i \sla{g} i$
for $g \in G$. Also let $\underline{U(1)}$ be the trivial `bundle of groups' on $X$ \cite{ref:moerdijk1}; as a groupoid it has
objects the elements $i \in X$ with hom-sets given by $\Hom(i,i) = U(1)$ and $\Hom(i,j) = \emptyset$.

\begin{figure}[t]
\center
\ig{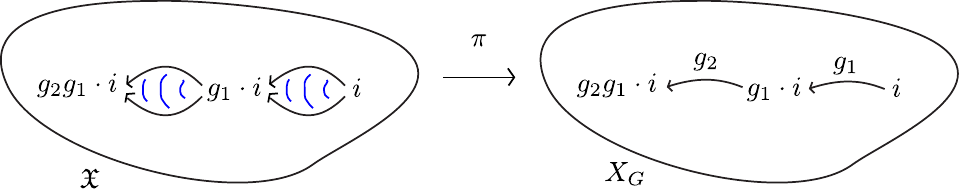}
\caption{\label{gerbefig}An equivariant gerbe.}
\end{figure}

\begin{defn} An {\em $G$-equivariant gerbe} is a central extension of the action groupoid of a discrete $G$-set $X$:
 \[
   \underline{U(1)} \stackrel{i}{\hookrightarrow} \X \sra{\pi} X_G.
 \]
A {\em metric} on the gerbe is an assignment of a positive real number $k_i$ to each object $i \in \X$, invariant under the action of $G$.
\end{defn}
By this we mean that $\X$ has the same objects as $X_G$, $\pi$ is a full surjective functor and $i$ is an
isomorphism onto the subgroupoid of arrows in $\X$ which project to identities in $X_G$ (see Figure \ref{gerbefig}).

We shall use a
non-calligraphic $X$ to refer to the underlying set of objects of an equivariant gerbe $\X$, and we say that the equivariant gerbe is {\em finite} if $X$ is a finite set.

\subsubsection*{Equivariant gerbes and cohomology}
We will write
 \[
  \X_{\sla{g} i } := \pi^{-1} (\sla{g} i)
   \]
for the $U(1)$-torsor of arrows in $\X$ which project to $\sla{g} i$ in $X_G$. A {\em section} of the gerbe is a
set-theoretic map
\[
  s \colon \text{Arr}(X_G)  \rightarrow  \text{Arr}(\X) \\
\]
such that $\pi \circ s = \id$. Choosing a section gives rise to a  $U(1)$-valued 2-cocycle $\phi \in Z^2(X_G,
U(1))$ on the groupoid $X_G$, in the sense of \cite{ref:simon}. One defines $\phi$ by
 \[
  s(\sla{g'} g \cdot x) \circ s(\sla{g} x) = \phi_x (g', g)
  s(\sla{g'g} x).
 \]
Choosing a different section $s'$ will change $\phi$ by a coboundary, so that an equivariant gerbe $\X$ gives
rise to a cohomology class $c_\X \in H^2(X_G, U(1))$.

\subsubsection*{Tensor product of equivariant gerbes} If $\X'$ and $\X$ are equivariant gerbes equipped with metrics, then their {\em
tensor product} $\X' \otimes \X$ is the equivariant gerbe with metric whose object set is the cartesian product $X' \times X$, whose $g$-graded morphisms are the tensor product of those of $\X'$ and $\X$,
 \[
  (\X' \otimes \X)_{\sla{g} (\mu, i)} := \X'_{\sla{g} \mu} \otimes \X_{\sla{g} i},
  \]
and whose metric is the product metric. Also, if $\X$ is an equivariant gerbe, we write $\overline{\X}$ for the equivariant gerbe having the same underlying groupoid as $\X$ but with the conjugate action of $U(1)$ on its hom-sets.

\subsection{Example \label{exeqgerbes}} Our main examples of equivariant gerbes are those arising from unitary 2-representations of $G$, but here is another example. Suppose $M$ is a
smooth connected manifold. Consider the groupoid $\mathcal{P}_M$ whose objects are $U(1)$-bundles with
connection $(P, \nabla)$ over $M$, and whose morphisms $f \colon (P, \nabla) \rightarrow (P', \nabla')$ are
diffeomorphisms $f \colon P \rightarrow P'$ which respect the action of $U(1)$ and which preserve the connection. If
there {\em is} an isomorphism from $(P, \nabla)$ to $(P', \nabla')$, then any other isomorphism must differ from
it by a constant factor in $U(1)$, since the maps must preserve parallel transport. Thus the hom-sets in
$\mathcal{P}_M$ are $U(1)$-torsors.

Now suppose a group $G$ acts on $M$ by diffeomorphisms. Let $\Pic^\nabla(M)$ denote the isomorphism classes in
$\mathcal{P}_M$ (it is given by a Deligne cohomology group), and suppose one chooses distinguished
representatives $(P, \nabla)_c$ for each isomorphism class $c \in \Pic^\nabla(M)$.  The group $G$ acts on
$\mathcal{P}_M$ by push-forward and hence on $\Pic^\nabla(M)$. This gives rise to
an associated equivariant gerbe $\X$ over $\Pic^\nabla(M)_G$ by the Grothendieck construction: the objects of
$\X$ are line bundles $c \in \Pic^\nabla(M)$ while the $g$-graded morphisms are given by
 \[
  \X_{\sla{g} i} := \Hom_{P_M} ( (L, \nabla)_{g \cdot c}, \,\, g_* (L, \nabla)_c).
 \]
Composition of $f_2 \in \X_{\sla{g_2} g_1 \cdot c}$ and $f_1 \in \X_{\sla{g_1} c}$ is defined by
 \[
 f_2 \diamond f_1 := \alpha_{g_2}(f_1) \, \circ \, f_2,
  \]
where we have (harmlessly) left out the canonical isomorphisms ${g_2}_* {g_1}_* (L, \nabla) \cong (g_2 g_1)_* (L, \nabla)$. This is simply a reformulation of the ideas in the second chapter of Brylinski \cite{ref:brylinski}. In
particular, it makes it clear that the automorphism group of any line bundle $c \in \X$ is a central extension of the
subgroup $H \subseteq G$ which fixes the isomorphism class $(P, \nabla)_c$. Many interesting central extensions
of groups arise in this way.

\subsection{Equivariant gerbes as a 2-category} In this subsection we define the 2-category of equivariant gerbes.

\subsubsection*{Unitary vector bundles over equivariant gerbes} A {\em unitary equivariant vector bundle} $E$ over an
equivariant gerbe $\X$ is a functor $E \colon \X \rightarrow \Hilb$ which maps morphisms in $\X$ to unitary maps in
$\Hilb$. A morphism $\theta \colon E \rightarrow E'$ of equivariant vector bundles over $\X$ is a natural
transformation; we write $\Hilb(\X)$ for the category of unitary equivariant vector bundles over $\X$.

If we choose a set-theoretic section $s$ of $\X$, giving rise to a 2-cocycle $\phi \in Z^2(X_G, U(1))$, then a
unitary equivariant vector bundle over $\X$ can be regarded as a $\phi$-twisted equivariant vector bundle
$\hat{E}$ over $X_G$ in the sense of \cite{ref:simon}, using the prescription $\hat{E}(\sla{g} i) = E(s(\sla{g} i))$. Functoriality of $E$ means
that $\hat{E}$ is $\phi$-twisted functorial,
 \[
  \hat{E}(\sla{g_2} g_1 \cdot i) \hat{E}(\sla{g_1} i) = \phi_i(g_2, g_1) \hat{E}(\la{g_2g_1} i).
 \]

\subsubsection*{The definition}
\begin{figure}[t]
\centering
\ig{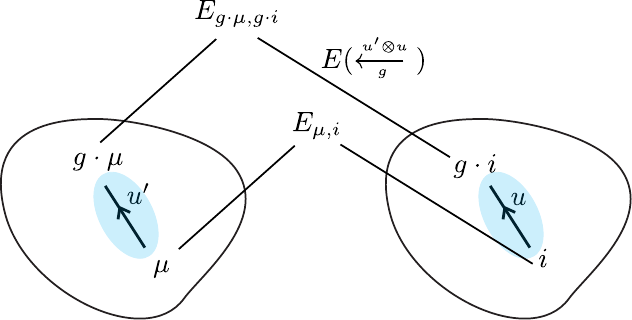}
\caption{\label{morfigg}A morphism of equivariant gerbes.}
\end{figure}
\begin{defn} For a finite group $G$, the 2-category $\Gerbes(G)$ of {\em finite $G$-equivariant gerbes} is defined as follows. An object is a
finite $G$-equivariant gerbe $\X$ equipped with a metric. The category of morphisms from $\X$ to $\X'$ is the category of unitary
equivariant vector bundles over $\X' \otimes \overline{\X}$,
\[
 \Hom(\X, \X') := \Hilb(\X' \otimes \overline{\X}).
  \]
\end{defn}
We refer the reader to Figure \ref{morfigg}. Composition of 1-morphisms works as follows --- the composite of $ \X'' \sla{E'}
\X' \sla{E} \X$ is the unitary equivariant vector bundle $E' \circ E$ over $\X'' \otimes \overline{\X}$ with fibers given by the weighted sum of Hilbert spaces
 \be \label{2kerncomp}
  (E' \circ E)_{\kappa, i} = \hat{\displaystyle \bigoplus_{\mu \in X'}} k_\mu \; E'_{\kappa, \mu} \otimes E_{\mu, i}.
 \ee
That is to say, the inner product is given on the homogenous components by
 \[
  (v_1 \otimes w_1, v_2 \otimes w_2) = k_\mu (v_1, v_2)(w_1, w_2)
 \]
where $k_\mu$ is the scale factor on $\mu \in X'$ (this formula should be compared with the formula for composition of kernels in \eqref{kerncomp}, and underscores the fact that the scale factors $k_\mu$ should be thought of as a metric). The equivariant maps on $E' \circ E$ are given by
 \be \label{cdeff}
  (E' \circ E)( \lag{u'' \otimes u}{g}) = \hat{\displaystyle \bigoplus_{\mu \in X'}} E'( \lag{u'' \otimes u'}{g}) \otimes E(\lag{u'
  \otimes u}{g}),
  \ee
where $u'' \in \X''_{\sla{g} \kappa}$, $u \in \X_{\sla{g} i}$ and the choices $u' \in \X'_{\sla{g} \mu}$ can be made
arbitrarily since the formula is invariant under $u' \mapsto e^{i \theta} u'$. The hats on these direct sums are there to indicate that they are not {\em formal} direct sums (which would require an ordering) but are rather defined via a geometric pull-push formula identical to the formula for composition of kernels in the derived category context (see \cite{ref:ganter_kapranov_rep_char_theory, ref:caldararu_willerton}). Horizontal and vertical composition of 2-morphisms works in a similar way. Note that $\Gerbes(G)$ is not a {\em strict} 2-category, but that will not concern us here (see \cite{ref:bartlett_details} for more details).

\subsubsection*{Classification of equivariant gerbes} We say two equivariant gerbes are {\em equivalent} if they are equivalent in the 2-category $\Gerbes(G)$; if they have metrics then we say they are {\em isometrically equivalent} if the support of the vector bundle $E\colon \X \rightarrow \X'$ furnishing the equivalence pairs together objects with the same scale factor. The following classification result is useful to bear in mind (see \cite{ref:bartlett_details} for more details).
\begin{prop} \label{equivgerbes} Suppose $\X$ and $\X'$ are equivariant gerbes equipped with metrics. The following are equivalent:
 \begin{enumerate}
  \item $\X$ is isometrically equivalent to $\X'$.
  \item There exists an isomorphism of $G$-sets $f \colon X \rightarrow X'$, preserving the scale factors, such that $c_{\X} = f^*(c_{\X'})$ as cohomology classes in $H^2(X_G, U(1))$.
 \end{enumerate}
\end{prop}

\subsection{U(1)-bundles and line bundles\label{uonebund}}
In this subsection we define $U(1)$-bundles and line bundles on groupoids and their spaces of sections.
\subsubsection*{Hermitian lines and $U(1)$-torsors} Recall that a {\em hermitian line} is a one-dimensional complex vector space with inner product, and a {\em $U(1)$-torsor} is a set equipped with a free and transitive action of $U(1)$. We write $\UOneTor$ for the category of $U(1)$-torsors and equivariant maps, and $\mathcal{L}$ for the category of hermitian lines and linear maps. To a $U(1)$-torsor $P$ we can associate a hermitian line $P_\mathbb{C}$ by taking the quotient of the cartesian product $P \times \mathbb{C}$ under the identifications $(e^{i \theta} p, \lambda) \sim (p, e^{i \theta} \lambda)$. We write the equivalence class of $(p, \lambda)$ as $p \otimes \lambda$, and the inner product on the line $P_\mathbb{C}$ is defined by $(p \otimes \lambda, p' \otimes \lambda') = \frac{p'}{p} \overline{\lambda}\lambda'$. Similarly, given a hermitian line $L$ we can associate a $U(1)$-torsor by taking the elements of unit norm.

\subsubsection*{$U(1)$-bundles and line bundles on groupoids}
We define a {\em $U(1)$-bundle with connection} over a finite groupoid $\mathcal{G}$ to be a functor $P \colon
\mathcal{G} \rightarrow \UOneTor$. Similarly a {\em hermitian line bundle with unitary connection} over
$\mathcal{G}$ is a functor $L \colon \mathcal{G} \rightarrow \mathcal{L}$, such that all the maps $L(\gamma)$ are
unitary, where $\gamma$ is an arrow in $\mathcal{G}$. We can use the conventions in the previous paragraph to
convert $U(1)$-bundles with connection into hermitian line bundles with unitary connection, and vice-versa.

\subsubsection*{$U(1)$-bundles and 1-cocycles}
A {\em trivialization} of a $U(1)$-bundle is a choice $\lambda_a \in P_a$ for each $a \in \mathcal{G}$. Choosing a
trivialization gives rise to a $U(1)$-valued cocycle $\alpha \in Z^1(\mathcal{G}, U(1))$ (in the sense of \cite{ref:simon}) whose value on a
morphism $\gamma$ in $\mathcal{G}$ is defined by the equation
 \[
  \lambda_{\text{target}(\gamma)} = \alpha(\gamma)
  P(\gamma)(\lambda_{\text{source}(\gamma)}).
 \]

\subsubsection*{Flat sections of line bundles}
A {\em flat section} of a line bundle $L \colon \mathcal{G} \rightarrow \mathcal{L}$ is a choice $s_a \in L_a$ for each $a
\in \mathcal{G}$, such that $s( \text{target}(\gamma)) = L(\gamma) s( \text{source}(\gamma))$ for all arrows
$\gamma \in \text{Arr} \, \mathcal{G}$. The space of flat sections of $L$ is denoted $\Gamma(L)$. If $s$ and
$s'$ are flat sections, then their fibrewise inner-product $(s, s')_x$ is a 0-form on $\mathcal{G}$, and hence can be integrated with respect to the natural measure on a groupoid (see \cite{ref:simon}), so that the space of sections $\Gamma(L)$ is endowed with an inner product via
 \[
  ( s, s') = \int_{x \in \mathcal{G}} (s, s')_x := \sum_{x \in \mathcal{G}} \frac{(s(x), s'(x))}{|x \sra{\gamma}|}.
 \]

\subsection{Transgression and twisted characters\label{tgess}} In this subsection we define the {\em transgressed line bundle} of an equivariant gerbe as a certain line bundle over the loop groupoid. Then we state the theorem of Willerton \cite{ref:simon} which identifies the space of isomorphism classes of equivariant vector bundles over the gerbe as the space of sections of this line bundle.

\subsubsection*{Line bundles from transgression of equivariant gerbes \label{transsec}}
Suppose $\X$ is an equivariant gerbe with underlying $G$-set $X$. Recall the notion of the {\em loop groupoid} $\Lambda X_G$ from Section \ref{defsec} --- the objects of the loop groupoid are `loops' in $X_G$ which we write as
\[
( \Fix{i}{x})
\]
and the morphisms are given by conjugation, which we write as
\[
(\Fixx{g\cdot i}{gxg^\mi}{5.7}{3}) \sla{g} (\Fix{i}{x}).
\]
We define the {\em transgressed $U(1)$-bundle} of $\X$ as the functor
 \begin{eqnarray*}
  \tau(\X) \colon \Lambda X_G & \longrightarrow & \UOneTor\\
  ( \Fix{i}{x}) & \mapsto & \X\sFix{i}{x} \\
  (\Fixx{g\cdot i}{gxg^\mi}{5.7}{3}) \sla{g} (\Fix{i}{x})  & \mapsto & u \mapsto vuv^\mi
  \end{eqnarray*}
where $v \in \X_{\sla{g} i}$ is an arbitrary choice; the formula is clearly independent of
this choice. The associated hermitian line bundle $\tau(\X)_\mathbb{C}$ is known as the {\em transgressed line bundle}.

\subsubsection*{Twisted characters of equivariant vector bundles\label{chsec}}
Suppose $E \colon \X \rightarrow \Hilb$ is a unitary equivariant vector bundle over an equivariant gerbe $\X$. The
{\em twisted character} (or just {\em character} for short) of $E$ is a flat section of the transgressed line bundle,
 \[
  \chi_E \in \Gamma_{\Lambda X_G} (\tau(\X)_\mathbb{C}).
 \]
It is defined by setting
 \[
  \chi_E ( \Fix{i}{x} ) = u \otimes \Tr E(\lag{u}{x})^*
 \]
where $u$ is any morphism $u \in \X \sFix{i}{x}$; the choice of $u$ doesn't matter since the formula is
invariant under $u \mapsto e^{i \theta} u$. We then have the following important theorem.
\begin{thm}[{Willerton \cite[Thm 11]{ref:simon}}] \label{Simthm} The twisted character map is a unitary isomorphism from the complexified Grothendieck group of isomorphism classes of unitary equivariant vector bundles over an equivariant gerbe $\X$ to the space of flat sections of the transgressed line bundle:
 \[
  \chi \colon [\Hilb_G (\X)]_\mathbb{C} \ra{\cong \, \, \, \,} \Gamma_{\Lambda X_G} (\tau(\X)_\mathbb{C}).
  \]
\end{thm}

\section{The geometric character of an equivariant gerbe\label{gch}}
This section is the geometric analogue of Section \ref{charsecc}: we define how to take the {\em geometric character} of a $G$-equivariant gerbe equipped with a metric in order to obtain a unitary equivariant vector bundle over $G$, and we show how to make this construction {\em functorial} with respect to morphisms of equivariant gerbes. We show that the geometric character descends to a functor from the homotopy category of $\Gerbes(G)$ to the category of equivariant vector bundles over $G$, and we use Theorem \ref{Simthm} to show that functor is unitarily fully faithful after one tensors the hom-sets with $\mathbb{C}$.

\subsection{The definition\label{deftransg}}
\begin{figure}[t]
\centering
\ig{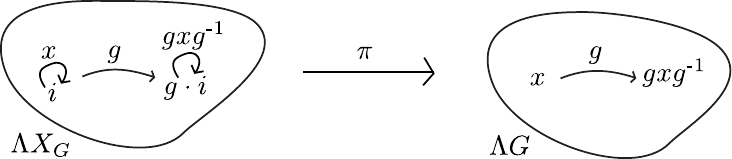}
\caption{\label{natfib}The fibration of loop groupoids associated to a $G$-set $X$.}
\end{figure}
For every $G$-set $X$ there is a natural fibration of loop groupoids
$\pi \colon \Lambda X_G \rightarrow \Lambda G$ whose fibers are the common fixed points (see Figure \ref{natfib}). If $X$ has a metric (an equivariant assignment of a positive real number $k_i$ to each $i \in X$) then we can
push-forward a unitary equivariant line bundle $L$ over $\Lambda X_G$ to a unitary equivariant vector bundle $\pi_* (L)$ over $G$ by taking the space of sections over the fixed points, as follows. The fiber of $\pi_*(L)$ at $x \in G$ is the Hilbert space
 \[
   \Sections (L|_{\FFix(x)}) = \hat{\displaystyle \bigoplus_{i \in \FFix{x}}} k_i \;  L \sFix{i}{x}
 \]
of sections of $L$ over the fixed points of $x$. That is, a vector $\psi \in \pi_* (L)_x$ is an assignment
  \[
   ( \Fix{i}{x}) \mapsto \psi (\Fix{i}{x}) \in L \sFix{i}{x}
  \]
where $i$ ranges over the fixed points of $x$, with the inner product on these sections is given by
 \[
  \langle \psi, \psi' \rangle = \sum_{i \in \FFix(x)} k_i \bigl( \psi(\Fix{i}{x}), \psi'(\Fix{i}{x}) \bigr).
 \]
As $x$ ranges over $G$, these vector spaces become a unitary equivariant vector bundle over $G$ via the natural $G$-action
 \[
 g \cdot \psi ( \Fixx{g\cdot i}{gxg^\mi}{5.7}{3}) = \psi(\Fix{i}{g}).
\]
In particular, if $\X$ is an equivariant gerbe equipped with a metric we write $\ch (\X) := \pi_* (\tau(\X)_\mathbb{C})$
for the push-forward of the transgressed line bundle of $\X$, and we call $\ch (\X)$ the {\em geometric character}
of $\X$.

\subsection{Functoriality for the geometric character\label{fgerbes}} We can make the geometric character functorial by
defining its action on morphisms $E \colon \X \rightarrow \X'$ of equivariant gerbes. Given such a morphism we can define for each $x\in G$
a linear map
 \begin{align*}
  \ch(E)_x \colon \Sections(\tau(\X)_\mathbb{C} |_{\FFix_X(x)}) & \rightarrow \Sections(\tau(\X')_\mathbb{C}
  |_{\FFix_{X'}(x)}) \\
   \psi & \mapsto \ch(E)_x (\psi)
  \end{align*}
by integrating the trace of $E$ over the fixed points of $x$. It is easiest to give this map in terms of its matrix elements between orthonormal bases of sections $\{ \psi_i \in \ch(\X)_x\}$ and $\{\psi'_\mu \in \ch(\X')_x\}$ consisting of sections which are localized over fixed points $i \in \FFix_X(x)$ and $\mu \in \FFix_{X'}(x)$ respectively:
\[
 \psi_i (\Fix{j}{x}) = \delta_{ij} \, u \otimes \frac{1}{\sqrt{k_i}} \in (\X \sFix{i}{x})_\mathbb{C}, \quad  \psi'_\mu (\Fix{\nu}{x}) = \delta_{\mu \nu} \, u' \otimes \frac{1}{\sqrt{k_\mu}} \in (\X' \! \ssFix{\mu}{x}{2}{0.5})_\mathbb{C}.
\]
%XXX discuss how 2-character distinguishes 2-reps
%XXX give some formulas for hom-sets using integration.
In terms of such a basis the matrix elements of $\ch(E)_x$ are defined as
 \be \label{matelements}
  \langle \psi'_\mu , \, \ch(E)_x \psi_i \rangle = \Tr E (\lag{u' \otimes u}{x})^*.
  \ee
Note that this definition of does not depend on the choices we made for $u' \in \X' \! \ssFix{\mu}{x}{2}{0.5}$ and $u \in \X \sFix{i}{x}$.

As before, we write $[\Gerbes(G)]$ for the {\em homotopy category} of $\Gerbes(G)$, and we write $[\Gerbes]_\mathbb{C}$ for the category whose hom-sets are the complexified Grothendieck groups of the hom-categories in $\Gerbes(G)$.
 \begin{thm}\label{geomchthm} The assignment
 \begin{align*}
  \chi \colon \Gerbes(G) & \quad \longrightarrow \quad \Hilb_G(G) \\
  \X & \quad \mapsto \quad \ch(\X) \\
  \X \sra{E} \X' & \quad \mapsto \quad \ch(\X) \ra{\ch(E)} \ch(\X')
 \end{align*}
is functorial with respect to composition in $\Gerbes(G)$, and only depends on the isomorphism class of $E$, and
thus descends to a {\em geometric character functor}
  \[
   \ch \colon [\Gerbes(G)] \rightarrow \Hilb_G(G).
   \]
Moreover, after tensoring the hom-sets in $[\Gerbes(G)]$ with $\mathbb{C}$ the associated functor
  \[
   \ch_\mathbb{C} \colon [\Gerbes(G)]_\mathbb{C} \rightarrow \Hilb_G(G)
  \]
 is unitarily fully faithful.
\end{thm}
\begin{proof}
Functoriality follows from the fact that trace is multiplicative on tensor products and `localizes' on the fixed
points. In other words, if $\X'' \la{E'} \X' \la{E} \X$ are morphisms of equivariant gerbes, and if $u \in \X
\sFix{i}{x}, u'' \in \X'' \sFix{\alpha}{x}$, then by the definition \eqref{cdeff} of composition in $\Gerbes(G)$,
\begin{align*}
 \Tr ( E' \circ E) (\lag{u'' \otimes u}{x}) &= \sum_{\mu \in X'} \Tr E' (\lag{u'' \otimes u'}{x}) \Tr E(\lag{u'
 \otimes u}{x}) \\
 &= \sum_{\mu \in \FFix_{X'}(x)} \Tr E'(\lag{u'' \otimes u'}{x}) \Tr E(\lag{u'
 \otimes u}{x}).
\end{align*}
Also the trace of an equivariant vector bundle only depends on isomorphism class; this gives the first
part of the proposition.

To prove the second part of the proposition, we show that the action of the geometric character on morphisms is
just a rearrangement of the twisted character map from Section \ref{chsec}, which is known to be unitary by Theorem \ref{Simthm}. Namely, we claim we have the following commutative diagram:
 \[
  \xymatrix{ \Hom_{[\Gerbes(G)]} (\X, \X') := [\Hilb (\X' \otimes \bar{\X})] \ar[r]^-{\chi} \ar[dr]_{\ch} &
  \Gamma_{\Lambda (X' \times X)_G} (\tau(\X' \otimes \bar{\X})_\mathbb{C}) \ar[d]^\wedge \\
  & \Hom_{\Hilb_G(G)} (\ch(\X), \ch(\X'))}
 \]
The rearrangement map $\hat{}$ (the downwards arrow) works as follows. A section $\xi$ of the transgressed line bundle $\tau(\X' \otimes
\bar{\X})_\mathbb{C}$ is something which assigns to every simultaneous fixed point an element of the appropriate hermitian line:
 \[
  \Fixx{(\mu, i)}{x}{4}{3} \quad \stackrel{\xi}{\mapsto} \quad u' \otimes u \otimes \lambda \; \in \; \X' \! \ssFix{\mu}{x}{2}{0.5} \otimes \overline{\X} \sFix{i}{x} \otimes \mathbb{C}.
 \]
Using the metrics on the gerbes, this can be regarded as a map of hermitian lines
\begin{align*}
 \hat{\xi} \colon (\X \sFix{i}{x})_\mathbb{C} \rightarrow (\X' \! \ssFix{\mu}{x}{2}{0.5})_\mathbb{C}
 \intertext{via the formula}
 \frac{1}{\sqrt{k_i}} u \otimes 1 \mapsto \frac{1}{\sqrt{k_\mu}} u' \otimes \lambda.
 \end{align*}
Using this correspondence, the section $\xi$ gives rise for each $x \in G$ to a linear map
 \[
  \hat{\xi}_x \colon \ch(\X)_x \rightarrow \ch(\X')_x.
 \]
The fact that $\xi$ was a {\em flat} section translates into the statement that the collection of maps $\hat{\xi}_x$ is {\em equivariant} with respect to the action of $G$. Moreover one can check that the map $\xi \mapsto \hat{\xi}$ is unitary with respect to the natural inner products involved (note that the twisted character map $\chi$ does not use the metrics on the gerbes, but $\ch$ and $\hat{}$ do), and also that the above diagram indeed commutes. Now we apply Theorem \ref{Simthm}, which says that after tensoring the left hand side with $\mathbb{C}$ the
character map $\chi$ is a unitary isomorphism. This gives the second statement of the proposition.
\end{proof}

\section{2-representations and equivariant gerbes\label{togsec}} In this section we show how to extract an equivariant gerbe from a marked unitary 2-representation, and similarly for morphisms and 2-morphisms, leading to a proof that the 2-category of unitary 2-representations is equivalent to the 2-category of equivariant gerbes. As we explained in the introduction, we encourage the reader to think of this as a `categorification' of the basic idea of geometric quantization (that representations of groups correspond to equivariant line bundles) in our simple discrete setting. Moreover, identifying unitary 2-representations with equivariant gerbes allows us to apply the integration technology of \cite{ref:simon}, resulting in some concrete formulas for the hom-sets in the complexified homotopy category of unitary 2-representations.

\subsection{Extracting equivariant gerbes from 2-representations\label{extract}} For technical reasons in this section we need to deal with {\em marked} 2-representations, by which we mean a unitary 2-representation on a 2-Hilbert space where a choice of representatives $e_i$ of the simple objects has been made. In any event, in practice many 2-Hilbert spaces arrive in this way; for instance it is common to have certain preferred choices for the irreducible representations of a group from the outset.

A unitary 2-representation $\alpha$ of $G$ on a marked 2-Hilbert space $H$ gives rise to an equivariant gerbe $\X$ by a variant of the {\em Grothendieck construction} (see for example \cite{ref:cegarra_et_al_graded_extensions}, or the original \cite{ref:grothendieck}), which we define as follows. The base set $X$ is the set of isomorphism classes of simple objects in
$H$; these are referred to as $i \equiv [e_i]$, etc. It inherits a $G$-action via
 \[
  g \cdot i := [\alpha_g (e_i)]
 \]
where $e_i$ is the distinguished simple object in the isomorphism class $[e_i]$. The equivariant gerbe $\X$ has the same objects as $X_G$, with the $G$-graded hom-sets given by
 \[
  \X_{ \sla{g} i} := \text{uIso}(e_{g \cdot i}, \alpha_g (e_i)).
  \]
Here ``uIso" refers to the unitary arrows in the 2-Hilbert space $H$ --- we warn the reader that such arrows do {\em not} lie on the unit circle in the hermitian line $\Hom(e_{g \cdot i}, \alpha_g(e_i))$, but rather on the circle with radius $\sqrt{k_i}$. The composite of $v$ and $v' \in \X_{\sla{g'} g \cdot i}$ is defined as
 \be \label{gerbecomp}
  v' \diamond v = \phi(g', g)_{e_i} \, \alpha_{g'}(v) \, v',
 \ee
and the identity morphisms $1_i \in \X_{\sla{e} i}$ are given by the unit isomorphisms $\phi_{e_i} (e)$. We have used the `$\diamond$' symbol above to stress that this is not {\em ordinary} composition of arrows. We will write the inverse of $v$ with respect to this composition law as $v^{\underline{\mi}}$.

\subsection{Examples\label{gexam}} We now extract the associated equivariant gerbes from some of the examples of unitary 2-representations we gave in Section \ref{exsec}.
\subsubsection{Automorphisms of groups} The equivariant gerbe $\X$ arising from the 2-representation of $G \subseteq{\Aut(K)}$ on $\Rep(K)$ works as follows. Firstly one chooses distinguished irreducible representations $V_i$ of $K$. The underlying $G$-set of $\X$ is the set of isomorphism classes of irreducible representations (which identifies noncanonically with the conjugacy classes of $G$), and the graded hom-sets are
 \[
  \X_{ \sla{g} i} = \uIso_K(V_{g \cdot i}, V_i^g).
 \]
In particular the $U(1)$-torsors above the fixed points are
 \[
  \X \sFix{i}{g} = \uIso_K(V, V^g) \subset U(V)
 \]
and there is no {\em apriori} preferred section of these torsors, unless $g$ is an inner automorphism. Let us summarize this discussion:
 \begin{quote}
 If $G$ acts as automorphisms of a group $K$, then each irreducible representation $\rho$ of $K$ will carry a projective representation of the subgroup $G_0 \subset G$ which fixes $\rho$.
 \end{quote}

\subsubsection{The metaplectic representation} In this way the equivariant gerbe arising from the 2-representation of $SL_2(\mathbb{R}) \subset \Aut(\text{Heis})$ on the category of representations of the Heisenberg group has a single object whose automorphism group is precisely the $\text{\em metaplectic}^c$ {\em group} $Mp^c(2)$ (the metaplectic group is the nontrivial double cover of $SL_2(\mathbb{R})$, see \cite[pg 2]{ref:robinson_rawnsley}). In other words, although the 2-representation is strict its corresponding equivariant gerbe and hence the associated projective representation are actually nontrivial.

\subsection{Morphisms of gerbes from morphisms of 2-representations} By a morphism of marked 2-representations we just mean an ordinary morphism of the underlying 2-representations which pays no attention to the distinguished simple objects, and similarly for the 2-morphisms; we write the 2-category of marked unitary 2-representations of $G$ as $\TRep_m(G)$.

A morphism $\sigma \colon \alpha \rightarrow \beta$ of marked 2-representations gives rise to a morphism of equivariant gerbes $\langle \sigma \rangle \colon \X_\alpha \rightarrow \X_\beta$ in the following way. The vector bundle $\langle \sigma \rangle$ over $\X_\beta \otimes \overline{\X}_\alpha$ is defined to
have fibers
 \[
  \langle \sigma \rangle_{\mu, i} := \langle \mu | \sigma | i \rangle  \stackrel{\text{shorthand}}{\equiv} \Hom(e_\mu, \sigma(e_i)).
 \]
For $v \in (\X_\beta)_{\sla{g} \mu}$ and $u \in (\overline{\X}_\alpha)_{\sla{g}
i}$, the unitary maps
 \[
 \langle \mu | \sigma | i \rangle (\lag{v \otimes u}{g}) \colon  \langle \mu | \sigma | i\rangle
  \rightarrow \langle g \cdot \mu| \sigma | g \cdot i \rangle
  \]
send
 \be \label{eqdefnmor}
 \ba \xymatrix{\sigma(e_i) \\ e_\mu \ar[u]^\lambda} \ea \quad
 \mapsto \quad\ba \xymatrix{ & \sigma \alpha_g (e_i)
 \ar[dr]^{\sigma(u)^*} \\
 \beta_g \sigma(e_i)\ar[ur]^{\sigma_{e_i}(g)} & & \sigma(e_{g \cdot i})
 \\
 \beta_g (e_\mu) \ar[u]^{\beta_g (\lambda)} & & e_{g \cdot \mu}
 \ar[ll]^{v}} \ea
 \ee
The fact that this construction defines a functor $\langle \sigma \rangle \colon \X_\beta \otimes
\overline{\X}_\alpha \rightarrow \Hilb$ follows from the coherence diagrams for the natural
isomorphisms $\sigma(g)$. The conjugate of $\X_\alpha$ must be used because of the $\sigma(u)^*$ term occurring
above.

\subsubsection*{2-Morphisms}
Similarly a 2-morphism $\theta \colon \sigma \rightarrow \rho$ between morphisms of marked 2-representations gives rise to a morphism of equivariant vector bundles $\langle \theta \rangle \colon \langle \sigma \rangle \rightarrow \langle \rho \rangle$ whose components are just given by postcomposition with the components of the natural transformation $\theta$:
 \begin{align*}
  \langle \theta \rangle_{\mu, i} \colon \langle \mu | \sigma | i \rangle & \rightarrow \langle \mu | \rho | i \rangle \\
   f & \mapsto \theta_{e_i} \circ f.
 \end{align*}

\subsection{Equivalence of 2-categories \label{eqref}} This allows us to identify the 2-category of unitary 2-representations with the 2-category of equivariant gerbes. This result is not really new, since related elements of it can be found in \cite[Cor 6.21]{ref:elgueta}, and
some similar ideas also appear in \cite[pg 17]{ref:barrett_mackaay}.  These references however do not use the language of equivariant gerbes. By using this language we believe our formulation expresses the {\em geometry} of the situation in a cleaner way, firstly because it shows how this result can be regarded as the `categorification' of the geometric correspondence between ordinary unitary representations and equivariant line bundles, and also because the language of equivariant gerbes is quite refined (see for instance Theorem \ref{Simthm}) and enables us to understand 2-representations and their 2-characters in a much better way. Moreover we essentially work directly with the 2-Hilbert spaces themselves and not some co-ordinatized skeleton of them, a strategy which is likely to be important in more advanced geometric situations.

\begin{thm} The map
 \begin{align*}
  \TRep_m(G) & \quad \longrightarrow \quad \Gerbes(G) \\
  \alpha & \quad \mapsto \quad \X_\alpha \\
  \alpha \sra{\sigma} \beta & \quad \mapsto \quad \X_\alpha \sra{\langle \sigma \rangle} \X_\beta \\
  \sigma \sra{\theta} \rho & \quad \mapsto \quad \langle\sigma \rangle \sra{\langle \theta \rangle} \langle \rho
  \rangle
 \end{align*}
is functorial, and an equivalence of 2-categories. Moreover for each pair $\alpha, \beta$ of marked unitary 2-representations the functor
  \[
   \Hom_{\TRep_m(G)} (\alpha, \beta) \rightarrow \Hom_{\Gerbes(G)} (\X_\alpha, \X_\beta)
  \]
given by the above prescription is a strong unitary equivalence of 2-Hilbert spaces.
\end{thm}
\begin{proof} The reason that this 2-functor is an equivalence is because unitary 2-representations, the morphisms between them and the 2-morphisms between those are determined by their behaviour on the simple objects. Since it is really a {\em weak} 2-functor the main thing to check is that the `compositor'
 \begin{align*}
 \langle \sigma' \rangle \circ \langle \sigma \rangle & \rightarrow \langle \kappa | \sigma' \circ \sigma | i \rangle \\
  g \otimes f & \mapsto \sigma'(f) \circ g
 \end{align*}
is equivariant with respect to the definition of the maps $ \langle \sigma' \circ \sigma)( \lag{u'' \otimes u}{g})$ from \eqref{cdeff}, which indeed turns out to be the case. Moreover the use of the scale factors $k_\mu$ in the definition of composition \eqref{2kerncomp} ensures that the compositor is a {\em unitary} isomorphism. That it is a strong unitary equivalence on the level of hom-categories follows from expanding out the definitions of the inner products on each side. See \cite{ref:bartlett_details} for details.
\end{proof}
Combining this result with the classification of equivariant gerbes from Proposition \ref{equivgerbes} allows us to rederive some known results about 2-representations, but specialized to the unitary setting. We say that a 2-representation is {\em irreducible} if its associated equivariant gerbe has only a single orbit.
\begin{cor}[{Compare \cite[Ex. 3.4]{ref:ostrik}, \cite[Thm 7.5]{ref:elgueta}, \cite[Prop
7.3]{ref:ganter_kapranov_rep_char_theory}}] Irreducible unitary 2-representations of $G$ are classified up
to strong unitary equivalence by triples $(X, k, [\phi])$, where $X$ is a transitive $G$-set up to isomorphism, $k$ is a positive real number, and $[\phi]$ is an equivariant cohomology class $[\phi] \in H^2( X_G, U(1))$.
\end{cor}
\noindent We write $1$ for the trivial 2-representation of $G$ on $\Hilb$ (it is the unit object for the monoidal 2-category structure of $\TRep(G)$, but we will not discuss this here).
\begin{cor}[{Compare \cite[Cor 6.21]{ref:elgueta}, \cite[Ex. 5.1]{ref:ganter_kapranov_rep_char_theory}}]  The endomorphism category of the unit object is monoidally equivalent to the category of representations of
$G$,
 \[
  \End_{\TRep(G)} (1) \simeq \Rep(G).
 \]
More generally, if $\alpha$ is any one-dimensional 2-representation of $G$, then there is a unitary equivalence
of 2-Hilbert spaces
 \[
  \Hom_{\TRep(G)} (1, \alpha) \simeq \Rep^\phi (G)
 \]
where $\phi \in Z^2(G, U(1))$ is the group 2-cocycle obtained from choosing a section of $\X_\alpha$.
\end{cor}
We would like to stress however that using the geometric language of equivariant gerbes allows us to {\em go further} than these results --- because it enables us to use the technology of \cite{ref:simon}, giving us a concrete understanding of {\em all} the hom-sets in $[\TRep(G)]_\mathbb{C}$.
 \begin{cor} The space of morphisms between unitary 2-representations in $[\TRep(G)]_\mathbb{C}$ identifies as the space of flat sections of the transgressed line bundle over the loop groupoid of the product of their associated $G$-sets:
  \[
   \Hom(\alpha, \beta)_{[\TRep(G)]_\mathbb{C}} \cong \Gamma_{\Lambda (X_\beta \times X_\alpha)_G} (\tau (\X_\beta \otimes \overline{\X_\alpha})_\mathbb{C}).
  \]
Thus their dimensions are given by
  \[
   \dim \Hom (\alpha, \beta) = \int_{\Lambda^2 (X_\beta \times X_\alpha)_G} \tau^2(\X_\beta \otimes \overline{\X_\alpha}).
  \]
In particular, the dimension of the space of endomorphisms of an object computes as
 \[
  \dim \End (\alpha) = \frac{1}{|G|} |\{(i, j, g, h) : i,j  \in X_\alpha,  g,h \in G, i,j \in \FFix(g) \cap \FFix(h), gh=hg\}|.
 \]
\end{cor}

\section{The 2-character and the geometric character\label{2char2rep}}
In this section we prove our main result in this paper --- that the 2-character of a unitary 2-representation corresponds naturally to the geometric character of its associated equivariant gerbe, and hence the 2-character is a unitarily fully faithful functor at the level of the complexified homotopy category.
\begin{thm} \label{transthm} The 2-character of a marked unitary 2-representation is unitarily naturally
isomorphic to the geometric character (i.e. the push-forward of the transgression) of the associated equivariant
gerbe:
 \[
 \xymatrix @1 @C=0.1in{[\TRep_m(G)] \ar[dr]_{\chi} \ar[rr]^{\sim} && [\Gerbes(G)] \ar[dl]^{\ch} \\ & \Hilb_G(G)}
 \]
\begin{textblock}{0.3}(0.245,-0.0515)
  $   \ba \ig{xy1.pdf} \ea $
\end{textblock}
\noindent That is, there are unitary isomorphisms $\gamma_\alpha \colon \chi_\alpha \stackrel{\cong}{\rightarrow} \ch(\X_\alpha)$,
natural in $\alpha$.
\end{thm}
Combining this with our knowledge of the geometric character functor from Theorem \ref{geomchthm} gives
\begin{cor} The complexified 2-character functor
 \[
  \chi_\mathbb{C} \colon [\TRep(G)]_\mathbb{C} \rightarrow \Hilb_G(G)
  \]
is a unitarily fully faithful functor from the complexified homotopy category of unitary representations of $G$ to the category of unitary equivariant vector bundles over $G$.
\end{cor}
Proving this theorem involves expanding out the abstract higher-categorical definitions for the 2-character and checking that they have the appropriate geometric behaviour. We do this in three steps --- firstly we define the isomorphisms $\gamma_\alpha$, then we show that they are indeed morphisms of equivariant vector bundles, and then we show that they are natural in $\alpha$.

\subsubsection*{Defining the isomorphisms}
Given a marked 2-representation $\alpha$ the isomorphism of
equivariant vector bundles $\gamma_\alpha \colon \chi_\alpha \rightarrow \ch(\X_\alpha)$ is easy enough to write
down. For $x \in G$, the fibers of the 2-character compute, by definition, as
 \begin{align*}
  \chi_\alpha (x) &= \Nat (\id, \alpha_x) \\
   & \cong \{ (\theta_{e_i} \colon e_i \rightarrow \alpha_x(e_i))_{i \in \FFix(x)} \} ,
 \end{align*}
while the fibers of the geometric character are
 \begin{eqnarray*}
   \ch (\X_\alpha)(x) &=& \Sections \bigl(\tau(\X_{\alpha})_\mathbb{C} |_{\FFix(x)} \bigr) \\
    &=& \{ (\vartheta_i \in \text{uHom}(e_i, \alpha_x (e_i)) \otimes \mathbb{C})_{i \in \FFix(x)} \} .
  \end{eqnarray*}
So the fibrewise identification between these two complex lines is component-wise just the identification between a hermitian line and the line associated to its circle of radius $\sqrt{k_i}$. Recalling our conventions about $U(1)$-torsors and hermitian lines from Section \ref{uonebund}, the isomorphisms $\gamma_\alpha$ are given at the level of unitary elements $u \in \uIso(e_i, \alpha_x(e_i))$ by
 \begin{align*}
  \gamma_\alpha \colon \chi_\alpha(x) & \rightarrow \ch (\X_\alpha)(x) \\
   u & \mapsto u \otimes \sqrt{k_i}.
  \end{align*}
Moreover, recalling our conventions about the inner products on $\Nat(\id, \alpha_x)$ and $\Sections \bigl(\tau(\X_\alpha)_\mathbb{C} |_{\FFix(x)} \bigr)$ from Sections \ref{hilb2} and \ref{deftransg} respectively, one sees that $\gamma_\alpha$ is indeed a unitary isomorphism.

\subsubsection*{Verifying that the isomorphisms are equivariant}
The following lemma verifies that these fibrewise identifications are equivariant with respect to the action of
$G$, which for the 2-character is given by the string diagram formula from Section \ref{defsec} and for the geometric
character by the transgression formula from Section \ref{transsec}. In other words, $\gamma_\alpha$ is indeed a morphism
in $\Hilb_G(G)$.

\begin{lem} \label{transglem} Let $\alpha$ be a marked 2-representation of $G$.
\begin{enumerate}
 \item The equivariant maps for the 2-character $\chi_\alpha$ compute as the
transgression, in the sense that
 \begin{align*}
  \chi_\alpha(gxg^\mi \sla{g} x) \colon \Nat(\id, \alpha_x) & \rightarrow \Nat(\id, \alpha_{gxg^\mi}) \\
  \intertext{evaluates as}
  \chi_\alpha (gxg^\mi \sla{g} x)(\theta)_{e_{g \cdot i}} &= v \diamond \theta_{e_i} \diamond v^{\underline{\mi}}
  \end{align*}
where $v \colon e_{g \cdot i} \rightarrow \alpha_g (e_i)$ is any unitary arrow in the underlying 2-Hilbert space
$H_\alpha$, and $\diamond$ is the twisted composition law from Section \ref{extract}.
 \item Therefore the following diagram commutes:
 \[
  \ba \xymatrix{\chi_\alpha (x) \ar[d]_{\gamma_\alpha(x)} \ar[rr]^{\chi_\alpha (\sla{g} x)} && \chi_\alpha (gxg^\mi) \ar[d]^{\gamma_\alpha
  (gxg^\mi)} \\
  \ch (\X_\alpha)(x) \ar[rr]_{\ch(\X_\alpha)(\sla{g} x)} && \ch (\X_\alpha)(gxg^\mi)} \ea\, .
 \]
\end{enumerate}
\end{lem}
\begin{proof} (i) We need to evaluate the string diagram formula for the map
$\chi_\alpha(\sla{g} x)$,
 \[
  \ba \ig{c3.pdf} \ea \mapsto \ba \ig{c4.pdf} \ea.
   \]
The right hand side computes as
\begin{eqnarray*}
 \bigl[\chi_\alpha (gxg^\mi \sla{g} x)(\theta)\bigr]_{e_{g \cdot i}} &:=&  \phi(gx, g^\mi)_{e_{g \cdot i}} \phi(g,x)_{\alpha_{g^\mi} (e_g \cdot i)}
  \alpha_g (\theta_{\alpha_{g^\mi}(e_{g \cdot i})}) \phi^*(g, g^\mi)_{e_{g \cdot i}} \phi(e)_{e_{g \cdot i}}  \\
   &\stackrel{\text{(a)}}{=}& \phi(gx, g^\mi)_{e_{g \cdot i}} \phi(g,x)_{\alpha_{g^\mi} (e_g \cdot i)}
  \alpha_g (\theta_{\alpha_{g^\mi}(e_{g \cdot i})}) \alpha_g(v^{\underline{\mi}}) v  \\
  &\stackrel{\text{(b)}}{=}& \phi(gx, g^\mi)_{e_{g \cdot i}} \phi(g,x)_{\alpha_{g^\mi} (e_g \cdot i)}
  \alpha_g(\alpha_x(v^{\underline{\mi}})) \alpha_g(\theta_{e_i})v \\
 &\stackrel{\text{(c)}}{=}& \phi(gx, g^\mi)_{e_{g \cdot i}} \alpha_{gx}(v^{\underline{\mi}}) \phi(g,x)_{e_i} \alpha_g(\theta_{e_i})  v \\
    &\stackrel{\text{(d)}}{=}& (v \diamond \theta_{e_i}) \diamond v^{\underline{\mi}}
 \end{eqnarray*}
where (a) uses the expression for $v^{\underline{\mi}}$ with respect to the twisted composition law $\diamond$ from \eqref{gerbecomp}, (b) uses the naturality of $\theta$, (c)
uses the naturality of $\phi(g,x)$, and (d) again uses the  composition law $\diamond$ from
\eqref{gerbecomp}.

(ii) This is just (i), expressed more formally.
\end{proof}

\subsubsection*{Verifying naturality}
It remains to show that the isomorphisms $\gamma_\alpha \colon \chi_\alpha \rightarrow \ch (\X_\alpha)$ are natural
with respect to morphisms of 2-representations $\sigma \colon \alpha \rightarrow \beta$. This amounts to computing
the string diagram formula \eqref{stdiagmor} for $\chi(\sigma)$ --- the action of the {\em 2-character} on
morphisms --- and observing that it corresponds to the formula \eqref{matelements} for the behaviour of the
geometric character on morphisms, which was defined in terms of the complex conjugate of the {\em
ordinary} character $\chi_{\langle \sigma \rangle}$ of the equivariant vector bundle $\langle \sigma \rangle$. In
a slogan, `the 2-character on morphisms is the ordinary character'.

\begin{lem} Let $\sigma \colon \alpha \rightarrow \beta$ be a morphism of marked 2-representations of $G$.
\begin{enumerate}
 \item The matrix elements of $\chi(\sigma)$ compute as the complex conjugate of the trace of the associated equivariant vector bundle $\langle \sigma \rangle$,
  \[
   \langle \theta_\mu, \chi(\sigma)_x \theta_i\rangle = \Tr \Bigl( \langle \mu | \sigma | i \rangle
   (\lag{\theta_\mu \otimes \theta_i}{x}) \Bigr)^*.
  \]
 \item Therefore the isomorphism $\gamma_\alpha \colon \chi_\alpha \sra{\cong} \ch(\X_\alpha)$ is natural in
 $\alpha$, i.e. the following diagram commutes:
  \[
   \xymatrix{\chi_\alpha(x) \ar[r]^{\chi(\sigma)_x} \ar[d]_{\gamma_\alpha} & \chi_\beta(x) \ar[d]^{\gamma_\beta} \\
             \ch(\X_\alpha)(x) \ar[r]_{\ch(\langle \sigma \rangle)_x} & \ch(\X_\beta)(x)}
  \]
  \end{enumerate}
\end{lem}
\begin{proof}
(i) We need to evaluate the string diagram formula for the map $\chi(\sigma)_x$,
\[
\ba \ig{y25.pdf} \ea \mapsto \ba \ig{y31.pdf} \ea .
\]
This gives
 \[
 \chi(\sigma)_x(\theta)_{e_\mu} = \beta_g(\epsilon_{e_\mu}) \circ \sigma(g)_{\sigma^*(e_\mu)} \circ \sigma(\theta_{\sigma^*(e_\mu)}) \circ n_{e_\mu}
 \]
where $\epsilon \colon \sigma \sigma^* \Rightarrow \id_{H_\beta}$ is the counit of the adjunction $\sigma^* \dashv \sigma$ and $n \colon \id_{H_\beta} \Rightarrow \sigma \sigma^*$ is the unit of the adjunction $\sigma \dashv \sigma^*$. Recall that $n$ is not arbitrary but is determined in terms of $\epsilon$ by the even-handed structure on $\THilb$, given at the level of adjunction isomorphisms by sending $\varphi \mapsto *\, \varphi^* *$. We can evaluate everything explicitly if we choose a {\em *-basis}
 \[
  \{ a_p \colon e_\mu \rightarrow \sigma(e_i)\}_{p=1}^{\dim \langle \mu | \sigma | i \rangle}
 \]
for each hom-space $\langle \mu | \sigma | i \rangle$; that is, a basis satisfying
 \[
  a_p^* \circ a_q = \delta_{pq} \id_{e_\mu} \quad \text{and} \quad \sum_{\mu, p} a_p \circ a_p^* = \id_{\sigma(e_i)}.
 \]
Substituting everything in, one computes that for a natural transformation $\theta_i \in \chi_\alpha (x)$ supported on a single fixed point $i \in \FFix_X (x)$ (that is, the components $(\theta_i)_{e_j}$ of $\theta_i$ over the marked simple objects are zero unless $i=j$), we have:
 \[
    \chi(\sigma)_x(\theta_i)_{e_\mu} = \frac{k_i}{k_\mu} \sum_p \beta_x (a_p^*) \circ \sigma(x)^*_{e_j} \circ \sigma(\theta_{e_j}) \circ a_p.
\]
Notice how the scale factors $k_i$ and $k_\mu$ have entered this description --- this underscores our point that this map {\em uses the even-handed structure in an intrinsic way}. We can identify this combination of terms as the complex conjugate of the trace of the associated equivariant vector bundle $\langle \sigma \rangle$ as follows. Fix an orthonormal basis $\{\theta_i \in \chi_\alpha (x)\}$ and $\{\theta_\mu \in \chi_\beta(x)\}$ of natural transformations supported exclusively over fixed points $i \in \FFix_X (x)$ and $\mu \in \FFix_{X'} (x)$. Recalling the relevant inner product from Section \ref{hilb2}, the matrix elements in this basis are thus
 \begin{align*}
  \langle \theta_\mu, \, \chi(\sigma)_x \theta_i \rangle &= k_\mu (\theta_\mu, \, \chi(\sigma)_x \theta_i) \\
   &= \sum_p k_i (\theta_{e_\mu}, \,  \beta_x(a_p^*)\sigma(x)_{e_i}^* \sigma(\theta_{e_i}) a_p) \\
   &= \frac{1}{k_\mu} \sum_p (a_p, \sigma(k_i \theta_{e_i}^*) \sigma(x)_{e_i} \beta_x(a_p) k_\mu \theta_{e_\mu})^* \\
   &= \Tr \Bigl( \langle \mu | \sigma | i \rangle
   (\lag{\theta_\mu \otimes \theta_i}{x}) \Bigr)^*,
   \end{align*}
where the last step uses the definition of the equivariant vector bundle $\langle \sigma \rangle$ from
\eqref{eqdefnmor} (we needed to use $k_i \theta_{e_i}$ because $\theta_i$ was an orthonormal basis vector, so that $(\theta_{e_i}, \theta_{e_i}) = \frac{1}{k_i}$, and similarly for $\theta_{e_\mu}$). We also used the fact that the trace of a linear endomorphism $A$ of the Hilbert space $\langle \mu | \sigma | i \rangle$ can be expressed in terms of a $*$-basis $\{a_p\}$ as
 \[
  \Tr(A) = \frac{1}{k_\mu} \sum_p (a_p, Aa_p)
 \]
since we must account for the fact that the basis vectors $a_p$ are not orthonormal:
 \[
  (a_p, a_q) = (\id_{e_\mu}, a_p^* a_q) = \delta_{pq} (\id_{e_\mu}, \id_{e_\mu}) = \delta_{pq} k_\mu.
 \]

(ii) This follows immediately from comparing the matrix elements of $\chi(\sigma)$ above to the matrix elements
of $\ch(\langle \sigma \rangle)$ given in \eqref{matelements}.
 \end{proof}

This completes the proof of Theorem \ref{transthm}.


\begin{thebibliography}{77} % start the bibliography
\small              % put the bibliography in a small font
\raggedright
%      --- A recommended citation format ---

\bibitem{ref:ambrose} W. Ambrose, Structure theorems for a special class of Banach algebras, Trans. Amer. Math. Soc. 57 (1945), 364-386.

\bibitem{ref:baez_dolan_HDA0} J. C. Baez and J. Dolan, Higher-dimensional algebra and topological
quantum field theory, Jour. Math. Phys. 36 (1995), 6073-6105.

\bibitem{ref:baez_2_hilbert_spaces} J. C. Baez, Higher-dimensional algebra II: 2-Hilbert
spaces, Adv. Math. 127 (1997), 125-189. Also available as \qalg{9609018}.

\bibitem{ref:baez_lauda_2-groups} J. C. Baez and A. D. Lauda, Higher-dimensional algebra V: 2-Groups, Theory and Applications of Categories Vol 12 (2004) No. 14, 423-491. Also available as \mathQA{0307200}.

\bibitem{ref:baez_higher_schreier_theory} J. C. Baez, This Week's Finds in Mathematical Physics Week 223, \href{http://math.ucr.edu/home/baez/week223.html}{\tt http://math.ucr.edu/home/baez/week223.html}.

\bibitem{ref:baez_baratin_freidel_wise} J. C. Baez, A. Baratin, L. Freidel and D. Wise, Representations of 2-Groups on Higher Hilbert Spaces, in preparation.

\bibitem{ref:barrett_mackaay} J. Barrett and M. Mackaay,  Categorical
representations of categorical groups, Theory and Applications of Categories Vol. 16 (2006) No. 20, 529-557. Also available as \mathCT{math.CT/0407463}.

\bibitem{ref:behrend_xu} K. Behrend and P. Xu, Differentiable Stacks and Gerbes, available as \Math{0605694}.

\bibitem{ref:bartlett_even_handed} B. Bartlett, Even-handed structures on 2-categories, in preparation.

\bibitem{ref:bartlett_ord} B. Bartlett, A categorical formulation of the equivalence between unitary representations of groups and equivariant line bundles, in preparation.

\bibitem{ref:bartlett_details} B. Bartlett, PhD thesis (University of Sheffield), in preparation.

\bibitem{ref:brylinski} J-L. Brylinski, {\em Loop Spaces, Characteristic Classes and Geometric Quantization},
Progress in Mathematics volume 107, Birkhauser-Boston 1993.

\bibitem{ref:brylinski_mclaughlin} J-L. Brylinski and D. McLaughlin, The geometry of degree-four characteristic classes and of line bundles on loop spaces I, Duke Math. J. 75 (1994) No. 3,
603-638.

\bibitem{ref:caldararu_willerton} A. Caldararu and S. Willerton, The Mukai Pairing, I: A categorical
approach. Available as \Math{0707.2052}.

\bibitem{ref:carter_segal_macdonald} G. Segal, {\em Lie Groups}, taken from R. Carter, G. Segal and I. Macdonald, {\em Lectures on Lie Groups and Lie Algebras},
London Mathematical Society, Student Texts 32 (1995).

\bibitem{ref:cegarra_et_al_graded_extensions} A.M. Cegarra, A.R. Garz\'{o}n and A.R.-Grandjean, Graded extensions of categories, J. Pure Appl. Algebra 154 (2000) 117-141. Also available as \href{http://www.ugr.es/~anillos/Preprints/total.ps}{\tt http://www.ugr.es/~anillos/Preprints/total.ps}.

\bibitem{ref:crane_yetter} L. Crane and D.N. Yetter, Measurable Categories and 2-Groups, Applied Categorical Structures 13 (2005) No. 5-6, 501-516. Also available as \mathQA{0305176}.

\bibitem{ref:deligne} P. Deligne, Action du groupe des tresses sur une cat\'{e}gorie, Invent. Math. 128 (1997), 159-175.

\bibitem{ref:dijkgraaf_witten} R. Dijkgraaf and E. Witten, Topological gauge theories and group cohomology, Comm. Math. Phys. 129 (1990) 60-72.

\bibitem{ref:elgueta} J. Elgeueta, Representation theory of 2-groups on Kapranov and Voevodsky's 2-category
2Vect, Adv. Math. 213 (2007) No. 1, 53-92. Also available as \Math{0408120v2}.

\bibitem{ref:fantechi} B.~Fantechi et al, {\em Fundamental Algebraic Geometry: Grothendieck's FGA Explained}, AMS Mathematical Surveys and Monographs Volume 123 (2006).

\bibitem{ref:freed} D. Freed, Higher algebraic structures and quantization, Comm. Math. Phys. 159 (1994) No. 2, 343-398.

\bibitem{ref:freed2} D. Freed, Quantum groups from path integrals. Available as \qalg{9501025}.

\bibitem{ref:freed_teleman_hopkins} D. Freed, M. Hopkins and C. Teleman, Loop Groups and Twisted K-Theory II. Available as \Math{0511232}.

\bibitem{ref:ganter_kapranov_rep_char_theory} N. Ganter, M. Kapranov, Representation and character theory in 2-categories, Adv. Math. 217 (2008) No. 5, 2268-2300. Also available as \mathKT{0602510}.

\bibitem{ref:grothendieck} A. Grothendieck, {\em Cat\'{e}gories fibr\'{e}es et descent (SGAI) expos\'{e}
VI}, Lecture Notes in Mathematics Vol. 224 (1971) 145-194, Springer Berlin.

\bibitem{ref:Gurski} N.~Gurski, {\em An algebraic theory of tricategories}, Phd Thesis, University of Chicago (2006). Available at \href{http://gauss.math.yale.edu/~mg622/tricats.pdf}{\tt http://gauss.math.yale.edu/~mg622/tricats.pdf}.

\bibitem{ref:kapranov_voevodsky} M. M. Kapranov and V. A. Voevodsky, 2-categories
and the Zamolodchikov tetradedra equations, Proceedings of
Symposia in Pure Mathematics 56 (1994) 177-259.

\bibitem{ref:kirwin} W. D. Kirwin, Coherent States in Geometric Quantization, J. Geom. Phys. 57 (2007) No. 2, 531-548. Also available as \Math{0502026}.

\bibitem{ref:lauda} A.~Lauda, Frobenius algebras and planar open string topological field theories. Available as \Math{0508349}.

\bibitem{ref:leinster_basic_bicategories} T. Leinster, Basic Bicategories. Available as \mathCT{9810017}.

\bibitem{ref:joyal_street} A. Joyal and R. Street, The geometry of tensor calculus I, Advances in Math. 88 (1991) 55-112.

\bibitem{ref:moerdijk1} I. Moerdijk, Introduction to the Language of Stacks and Gerbes. Available as \mathAT{0212266}.

\bibitem{ref:morton} J. Morton, Extended TQFT's and quantum gravity, PhD thesis University of California Riverside (2007). Available as \arXiv{0710.0032}.

\bibitem{ref:murray_singer} M. Murray and M. Singer, Gerbes,  Clifford Modules and the Index Theorem, Annals of Global Analysis and Geometry 26 (2004) No. 4, 355-367.  Also available as \mathDG{0302096}.

\bibitem{ref:ostrik} V. Ostrik, Module categories, weak Hopf algebras and modular invariants, Transformation Groups 8 (2003) No. 2, 177-206. Also available as \Math{0111139}.

\bibitem{ref:robinson_rawnsley} P.L. Robinson and J.H. Rawnsley, {\em The metaplectic representation, $Mp^\mathbb{C}$ structures and geometric quantization}, Memoirs of the American Mathematical Society Vol 81 No. 410.

\bibitem{ref:street_string_diagrams} R. Street, Categorical structures, in {\em Handbook
of Algebra}, Vol. 1, ed. M. Hazewinkel, Elsevier Amsterdam (1995), 529-577.

\bibitem{ref:tu_xu_laurent_gengoux} J-L. Tu, P. Xu and C.
Laurent-Gengoux, Twisted K-theory of Differentiable Stacks, Annales Scientifiques de l'\'{E}cole Normale Sup\'{e}rieure 37 (2004) No. 6, 841-910. Also available as \Math{0306138v2}.

\bibitem{ref:simon} S. Willerton, The twisted Drinfeld double of a finite group via gerbes and
finite groupoids. Available as \mathQA{0503266}.

\bibitem{ref:spera} M. Spera, On K\"{a}hlerian coherent states, Proceedings of the International Conference on
Geometry, Integrability and Quantization, Varna Bulgaria 1999 (eds. I. Mladenov and G. Naber), Coral Press (2000), 241-256. Available at \href{http://www.bio21.bas.bg/proceedings/Proceedings_files/vol1content.htm}{\tt http://www.bio21.bas.bg/proceedings/Proceedings\_files/vol1content.htm}.

\bibitem{ref:woodhouse} N. M. J. Woodhouse, {\em Geometric Quantization}, Oxford Mathematical Monographs (1992).
\end{thebibliography}
\end{document}